\documentclass[12pt]{amsart}
\usepackage{amssymb}
\numberwithin{equation}{section}
\theoremstyle{plain}
\newcommand{\edd}{\end{document}}
\newtheorem{lemma}{Lemma}[section]
\newtheorem{corollary}[lemma]{Corollary}
\newtheorem{proposition}[lemma]{Proposition}
\newtheorem{theorem}[lemma]{Theorem}
\newtheorem{sublemma}[lemma]{Sublemma}
\theoremstyle{remark}
\newtheorem{remark}[lemma]{Remark}
\newtheorem{definition}[lemma]{Definition}
\newtheorem{conjecture}[lemma]{Conjecture}

\def\on{\operatorname}

\def\Lemma{\begin{lemma}}
\def\enlemma{\end{lemma}}
\def\Definition{\begin{definition}}
\def\Def{\begin{definition}}
\def\endefinition{\end{definition}}
\def\endef{\end{definition}}
\def\Cor{\begin{corollary}}
\def\encor{\end{corollary}}
\def\Sublemma{\begin{sublemma}}
\def\ensublemma{\end{sublemma}}
\def\Prop{\begin{proposition}}
\def\enproposition{\end{proposition}}
\def\enprop{\end{proposition}}
\def\Theorem{\begin{theorem}}
\def\entheorem{\end{theorem}}
\def\Remark{\begin{remark}}
\def\enrem{\end{remark}}
\def\Conj{\begin{conjecture}}
\def\enconj{\end{conjecture}}

\newenvironment{tenumerate}{
  \begin{enumerate}
  
  }{\end{enumerate}}

\newenvironment{fnumerate}{
  \begin{enumerate}
  
  }{\end{enumerate}}

\newenvironment{anumerate}{
  \begin{enumerate}
  
  }{\end{enumerate}}

\newcommand{\bnum}{\begin{tenumerate}}
\newcommand{\enum}{\end{tenumerate}}
\newcommand{\banum}{\begin{anumerate}}
\newcommand{\eanum}{\end{anumerate}}
\newcommand{\bfnum}{\begin{fnumerate}}
\newcommand{\efnum}{\end{fnumerate}}

\def\b{{\mathfrak{b}}}
\def\g{{\mathfrak{g}}}

\def\t{{\mathfrak{t}}}
\def\GL{{\mathfrak{L}}}

\def\C{{\mathbb{C}}}

\def\BQ{{\mathbb{Q}}}

\def\BZ{{\mathbb{Z}}}
\def\Z{{\mathbb{Z}}}

\def\O{{\mathcal O}}

\def\lam{\lambda}
\def\lan{\langle}
\def\ran{\rangle}
\def\nn{\nonumber}

\def\cl{\operatorname{{cl}}}
\def\mod{\operatorname{{mod}}}

\def\gr{^{\,\circ}}

\def\Ker{\mathop{\rm Ker\hskip.5pt}\nolimits}

\def\max{{\mathop{\rm{max}}}}

\def\re{{\rm re}}
\def\sHom{\mathop{\mathscr{H}\kern-3pt\hbox{\it om}\,}\nolimits}

\def\wt{{\on{wt}}}
\def\Wt{{\rm Wt}}
\newcommand{\var}{\on{var}}
\newcommand{\vphi}{\varphi}
\newcommand{\vpi}{\varpi}

\newcommand{\U}{{U_q(\g)}}
\newcommand{\iso}{%
\raise-.7pt\hbox{$\mathrel{\longrightarrow{\kern-14.5pt\raise3.5pt%
\hbox{$\sim$}}}$}\enspace}
\newcommand{\sqtimes}{\mathrel{\square\kern-7.6pt\raise1pt\hbox{$\times$}}}
\newcommand{\ssqtimes}{\mathrel{\square\kern-7pt\raise.05pt\hbox{$\times$}}}
\def\endeq{\end{eqnarray}}
\def\eneq{\end{eqnarray}}
\def\endeqn{\end{eqnarray*}}
\def\eneqn{\end{eqnarray*}}
\def\eq{\begin{eqnarray}}

\def\eqn{\begin{eqnarray*}}
\newcommand{\isoto}{\mathrel{\longrightarrow{\kern-18pt\raise3.5pt\hbox{$\sim$
}}}}

\def\Q{{\mathbb Q}}
\def\pf{\proof}

\def\tf{{\tilde f}}

\newcommand{\ol}{\overline}
\newcommand{\eps}{\varepsilon}
\newcommand{\vp}{\varphi}
\newcommand{\te}{{\tilde e}}
\newcommand{\TF}{\widetilde{F_i}{}^{(n)}}
\newcommand{\qbin}[3][i]{{\genfrac{[}{]}{0pt}{}{#2}{#3}}_{#1}}

\def\hb{\hfill\break}

\def\qbox#1{\hbox{\quad#1\quad}}

\def\Gg{{\mathfrak{g}}}

\def\Gt{{\mathfrak{t}}}

\def\BQ{{\mathbb Q}}

\def\BZ{{\mathbb Z}}
\def\Z{\BZ}
\def\lam{\lambda}
\def\Lam{\Lambda}
\def\vphi{\varphi}
\def\vp{\varphi}
\def\lan{\langle}
\def\ran{\rangle}

\def\ol{\overline}
\def\on{\operatorname}
\def\gge{>\kern-3pt>}

\def\cl{{\rm{cl}}}
\def\mod{{\rm{mod}}}

\def\Ker{\mathop{\rm Ker\hskip.5pt}\nolimits}

\def\re{{\rm re}}

\def\max{{\mathop{\rm{max}}}}
\def\varprojlim{\mathop{\vtop{\ialign{$##$\cr
\hfil{\fam0 lim}\hfil\cr\noalign{\nointerlineskip}%
{\leftarrow}\mkern-6mu\cleaders\hbox{$\mkern-2mu{-}\mkern-2mu$}\hfill
\mkern-6mu{-}\cr\noalign{\nointerlineskip\kern-.2326ex}\cr}}}}

\def\te{\tilde e}
\def\tf{\tilde f}

\def\aff{{\rm aff}}

\def\beq{\begin{eqnarray}}
\def\beqn{\begin{eqnarray*}}
\def\endeq{\end{eqnarray}}
\def\endeqn{\end{eqnarray*}}
\def\eq{\begin{eqnarray}}
\def\eqn{\begin{eqnarray*}}
\def\nn{\nonumber}
\def\proof{\noindent{\it Proof.}\quad}
\def\qed{\hspace*{\fill}{Q.E.D.}\par\medskip}

\def\U{U_q(\Gg)}
\def\Us{U'_q(\Gg)}
\def\Uf{U_q^-(\Gg)}
\def\Ue{U_q^+(\Gg)}

\def\rank{{\on{rank}}}

\def\clwt{\Gt^{*0}_\cl}

\def\semi{\hbox{$\,{\vrule height5.9pt depth.6pt}\kern-1.4pt\times$}}
\def\real{{\rm{re}}}
\def\tU{{\tilde U}_q(\Gg)}
\def\O{{\on{O}}}
\def\eps{\varepsilon}

\def\tD{\tilde\Delta}
\def\GL{{\on{GL}}}
\def\xio{\xi{\raise-2.7pt\hbox{${}_{0}$}}}
\def\refl#1{s{\raise-2.5pt\hbox{${}_{#1}$}}}

\def\hookdownarrow%
{{\raise3pt\hbox{$\cap$}\kern-9.9pt\raise-4pt\hbox{$\downarrow$}}}
\def\doublevrule
{\raise-5pt\hbox{$\vrule height15pt\hskip2pt\vrule height15pt$}}

\newcommand{\tQ}{{\widetilde Q}}
\newcommand{\Db}{\ol{\Delta}}

\newcommand{\ct}{{\widehat\otimes}}
\newcommand{\rct}{{\widetilde\otimes}}
\newcommand{\Up}{U_q^+(\Gg)}
\newcommand{\Um}{U_q^-(\Gg)}
\newcommand{\CT}{\Up\ct\Um}
\newcommand{\bt}{\ol{\otimes}}
\newcommand{\cbt}{\widehat{\ol{\otimes}}}
\newcommand{\good}{good}
\newcommand{\Good}{Good\ }
\newcommand{\Mh}{\widehat{M}}
\newcommand{\bi}{\begin{itemize}}
\newcommand{\ei}{\end{itemize}}
\newcommand{\hs}{\hspace*}
\newcommand{\vs}{\vspace}
\newcommand{\norm}{{\on{norm}}}
\newcommand{\univ}{{\on{univ}}}
\newcommand{\Rn}{R^\norm}
\newcommand{\Ru}{R^\univ}
\newcommand{\Rc}{R^{{\on{comb}}}}

\newcommand{\bio}{\bar\iota}
\newcommand{\cj}[1][]{{\on{c}_{#1}}}
\newcommand{\cju}{{\on{c}}^\univ}
\newcommand{\cjn}{{\on{c}}^\norm}
\newcommand{\qs}{{q_s}}

\newcommand{\ba}{\begin{array}}
\newcommand{\ea}{\end{array}}
\newcommand{\basic}{basic}
\newcommand{\dz}{{0^\vee}}
\newcommand{\bB}{\ol{B}}
\newcommand{\al}{\alpha}
\newcommand{\la}{\lambda}

\newcommand{\seps}{\varepsilon^*}
\newcommand{\sphi}{\varphi^*}
\newcommand{\set}{{\,;\,}}
\newcommand{\ccdot}{\,\cdot\,}
\newcommand{\mult}{\begin{multline}}
\newcommand{\enmult}{\end{multline}}

\title[Quantized
Affine Algebras]%
{On Level Zero Representations of Quantized
Affine Algebras}
\author[M.~Kashiwara]{Masaki Kashiwara}
\address[M.~Kashiwara]{Research Institute for Mathematical Sciences, Kyoto
University, Kyoto, 606--8502, Japan, L'Ecole Normale Sup\'erieure,
45 rue d'Ulm, 75230 Paris cedex 05, France
}
\thanks{This work benefits from
a ``Chaire Internationale de Recherche Blaise Pascal de l'Etat et de la
R\'egion d'Ile-de-France, g\'er\'ee par la
Fondation de l'Ecole Normale Sup\'erieure''.}

\keywords{crystal base, quantized affine algebra}
\subjclass{17B37}
\begin{document}
\begin{abstract}
We study the properties of level zero
modules over quantized affine algebras.
The proof of the conjecture
on the cyclicity of tensor products by Akasaka and the present author
is given.
Several properties of modules generated by extremal vectors
are proved.
The weights
of a module generated by an extremal vector 
are contained in the convex hull of the Weyl group orbit of 
the extremal weight.
The universal extremal weight module with level zero fundamental weight 
as an extremal weight is irreducible, and
isomorphic to the affinization
of an irreducible finite-dimensional module.
\end{abstract}
\maketitle




\tableofcontents

\section{Introduction}
In this paper, we study the level zero representations
of quantum affine algebras.
This paper is divided into three parts,
on extremal weight modules, on the conjecture
in \cite{AK} on the cyclicity of the tensor products of 
fundamental representations, and on the global basis of the Fock space.

 In \cite{modified},
as a generalization of highest weight vectors, 
the notion of extremal weight vectors
is introduced, and
it is shown that the universal module generated by
an extremal weight vector has favorable properties:
this has a crystal base, a global basis, etc.
The main purpose of the first part (\S~\ref{sec:rev1}---\S~\ref{sec:ext})
is to study such  modules in the {\em affine case}
and to prove the following two properties.
\bi
\item[(a)]
If a module is generated by an extremal vector with weight $\lam$,
then all the weights of this module are contained in the convex hull of
the Weyl group orbit of $\lam$.
\item[(b)]
Any module generated by an extremal vector
with a level zero fundamental weight 
$\varpi_i$ 
is irreducible, and
isomorphic to the affinization
of an irreducible finite-dimensional module $W(\varpi_i)$
(see Theorem~\ref{th:fundrep} and Proposition
~\ref{prop:fundrep} for an exact statement).
\ei

In the second part, we shall prove the following theorem
\footnote{M. Varagnolo--E. Vasserot
(Standard modules of quantum affine algebras, math.QA/0006084)
prove the same conjecture in the simply-laced case
by a different method.},
which is conjectured in \cite{AK}
and proved in the case of $A^{(1)}_n$ and $C^{(1)}_n$.

\vs{3pt}
\noindent
{\bf Theorem.}\quad{\em
If $a_\nu/a_{\nu+1}$ has no pole at $q=0$ {\rm ($\nu=1,\ldots,m-1$)}, then
$W(\varpi_{i_1})_{a_1}\otimes\cdots\otimes W(\varpi_{i_m})_{a_m}$
is generated by the tensor product of the extremal vectors.
}

\vs{3pt}
In the course of the proof, one uses the global basis on
the tensor products of the affinizations of $W(\varpi_{i_\nu})$,
especially the fact that the transformation matrix between
the global basis of the tensor products 
and the tensor products of global bases
is triangular.

Among the consequences of this theorem (see \S~\ref{sec:main}),
we mention here the following one.
Under the conditions of the theorem above,
there is a unique homomorphism up to a constant multiple
$$
W(\varpi_{i_1})_{a_1}\otimes\cdots\otimes
W(\varpi_{i_m})_{a_m}
\longrightarrow
W(\varpi_{i_m})_{a_m}\otimes \cdots\otimes
W(\varpi_{i_1})_{a_1},
$$
and its image is an irreducible $\Us$-module.
This phenomenon is analogous to
the morphism from the Verma module to the dual Verma module.
Conversely, combining with a result of Drinfeld (\cite{D}),
any irreducible integrable $\Us$-module is isomorphic
to the image for some
$\{(i_1,a_1),\ldots,(i_m,a_m)\}$.
Moreover, $\{(i_1,a_1),\ldots,(i_m,a_m)\}$
is unique up to a permutation.

In the third part (\S~\ref{sec:Fock}),
we prove the existence of the global basis on the Fock space.

\vs{5pt}
The plan of the paper is as follows.
In \S~\ref{sec:rev1}
--\S~\ref{sec:aff}, we review some of the known results of crystal bases.
Then, in \S~\ref{sec:ext}, we give a proof of
(a) and (b).

In \S~\ref{sec:exis}, we prove a sufficient condition
for a module to admit a global basis:
very roughly speaking, it is enough to have a global basis in the extremal
weight spaces.
In \S~\ref{sec:Rmatrix}, we review the universal $R$-matrix and
the universal conjugation operator.
After introducing the notion of \good\ modules
(rudely speaking, a module with a global basis),
we shall prove in \S~\ref{sec:main} 
the above theorem in the framework of
\good\ modules 

After preparations in \S~\ref{sec:comb}--\S~\ref{sec:energy}
on the combinatorial $R$-matrix and the energy function,
we shall prove in \S~\ref{sec:Fock} the properties of \good\ modules
which are postulated for the existence of the wedge products and 
the Fock space in \cite{Fock}.
Finally, we shall show that
the Fock space admits a global basis.
In the case of the vector representation of $\g=A^{(1)}_n$,
the global basis of the corresponding Fock space
is already constructed by
B. Leclerc and J.-Y. Thibon \cite{LT} (see also
\cite{LT2,VV}).

In the last section, we present conjectures on
the structure of $V(\lam)$.

\vs{8pt}
{\em Acknowledgements.}\quad
The author would like to thank Anne Schilling
who kindly provided a proof of the formula \eqref{eq:anne}.

\section{Review on crystal bases}\label{sec:rev1}

In this section, we shall review very briefly 
the quantized universal enveloping algebra and crystal bases.
We refer the reader to \cite{banff,quant,modified}.

\subsection{Quantized universal enveloping algebras}
We shall define the quantized universal enveloping algebra
$U_q(\g)$.
Assume that we are given the following data.
\begin{eqnarray*}
&&P:\text{a free $\BZ$-module (called a weight lattice)}\\
&&I:\text{an index set (for simple roots)}\\
&&\al_i\in P\ \text{for}\ i\in I\ \text{(called a simple root)}\\
&&h_i\in P^*=\text{Hom}_\BZ(P,\BZ)\ \text{(called a simple coroot)}\\
&&(\,\cdot\,,\,\cdot\,)\colon P\times P\to \BQ\ \text{a bilinear symmetric form.}
\end{eqnarray*}
We shall denote by $\langle \,\cdot\,,\,\cdot\,\rangle 
\colon  P^*\times P\to \BZ$ the
canonical pairing. 

The data above are assumed to satisfy the following axioms.
\begin{eqnarray}
&&(\al_i,\al_i)>0\quad\text{for any $i\in I$,}\\
&&\langle h_i,\la\rangle =\frac{2(\al_i,\lam)}{(\al_i,\al_i)}
\quad\text{for any $i\in I$ and $\lam\in P$,}\\
&&(\al_i,\al_j)\le 0\quad 
\text{for any $i$, $j\in I$ with $i\not= j$.}
\end{eqnarray}

Let us choose a positive integer  $d$
such that $(\al_i,\al_i)/2\in\Z\,d^{-1}$ for any $i\in I$.
Now let $q$ be an indeterminate and set 
\eq
&&\text{$K=\Q(\qs)$ where
$\qs=q^{1/d}$.}
\eneq

\begin{definition}\label{U_q(g)}
The quantized universal enveloping algebra $U_q(\g)$ is the algebra over
$K$ generated by the symbols $e_i,f_i\ (i\in I)$ and $q(h)\ (h\in d^{-1}P^*)$ 
with the following defining relations.
\begin{enumerate}
\item $q(h)=1$ for $h=0$.

\item $q(h_1)q(h_2)=q(h_1+h_2)$ for $h_1,h_2\in d^{-1}P^*$.

\item
$q(h)e_i\,q(h)^{-1}=q^{\langle h,\al_i\rangle}\,e_i\ $
and $\ q(h)f_i\,q(h)^{-1}=q^{-\langle h,\al_i\rangle}f_i\ $
for any $i\in I$ and $h\in d^{-1}P^*$.

\item\label{even} $\lbrack e_i,f_j\rbrack
=\delta_{ij}\dfrac{t_i-t_i^{-1}}{q_i-q_i^{-1}}$
for $i$, $j\in I$. Here $q_i=q^{(\al_i,\al_i)/2}$ and
$t_i=q(\frac{(\al_i,\al_i)}{2}h_i)$.

\item (Serre relation) For $i\not= j$,
\begin{eqnarray*}
&&\sum^b_{k=0}(-1)^ke^{(k)}_ie_je^{(b-k)}_i=\sum^b_{k=0}(-1)^kf^{(k)}_i
f_jf_i^{(b-k)}=0.
\end{eqnarray*}
Here $b=1-\langle h_i,\al_j\rangle$ and
\begin{eqnarray*}
e^{(k)}_i=e^k_i/\lbrack k\rbrack_i!\ ,&& f^{(k)}_i=f^k_i/\lbrack k\rbrack_i!\ ,\\
\lbrack k\rbrack_i=(q^k_i-q^{-k}_i)/(q_i-q^{-1}_i)\ ,
&&\lbrack k\rbrack_i!=\lbrack 1\rbrack_i\cdots \lbrack k\rbrack_i\,.
\end{eqnarray*}
\end{enumerate}
\end{definition}

Sometimes we need an algebraically closed field containing $K$, for
example
\eq \widehat K=\bigcup_n \C((q^{1/n})),
\eneq
and to consider $\U$ as an algebra over $\widehat K$.

We denote by
$\U_\Q$ the subalgebra of $\U$ over $\Q[\qs^{\pm1}]$ generated by
the $e_i^{(n)}$'s, the $f_i^{(n)}$'s ($i\in I$)
and $q^h$ ($h\in d^{-1}P^*$).

Let us denote by $W$ the Weyl group, the subgroup of $GL(P)$ generated
by the simple reflections $s_i$:
$s_i(\lam)=\lam-\lan h_i,\lam\ran\alpha_i$.

Let $\Delta\subset Q=\sum_i\Z\alpha_i$ be the set of roots.
Let $\Delta^\pm=\Delta\cap Q_\pm$ be the set of positive 
and negative roots, respectively.
Here
$Q_\pm=\pm\sum_i\Z_{\ge0}\alpha_i$.
Let $\Delta^\re$ be the set of real roots.
$\Delta_\pm^\re=\Delta_\pm\cap\Delta^\re$.

\subsection{Crystals}
We shall not review the notion of crystals,
but refer the reader to
\cite{banff,quant,modified}.
We say that a crystal $B$ over $\U$
is {\it a regular crystal} 
if, for any $J{\subset}I$
such that $\{\al_i\set i\in J\}$
is of finite-dimensional type,
$B$ is, as a crystal over $U_q(\Gg_J)$, 
isomorphic to the crystal bases associated with
an integrable $U_q(\Gg_J)$-module.
Here $U_q(\Gg_J)$ is the subalgebra of $\U$ generated by
$e_j$, $f_j$ ($j\in J$) and $q^h$ ($h\in d^{-1}P^*$).
By \cite{modified},
the Weyl group $W$ acts on any regular crystal.
This action $S$ is given by
\eqn
&&S_{s_i}b=
\begin{cases}
\tf_i^{\lan h_i,\wt(b)\ran}b
&\mbox{if $\lan h_i,\wt(b)\ran\ge 0$,}\\
\te_i^{-\lan h_i,\wt(b)\ran}b
&\mbox{if $\lan h_i,\wt(b)\ran\le 0$.}
\end{cases}
\endeqn

Let us denote by
$\Uf$ (resp. $\Ue$) 
the subalgebra of $\U$ generated by the $f_i$'s
(resp. by the $e_i$'s).
Then $\Uf$ has a crystal base denoted by
$B(\infty)$ (\cite{quant}). The unique weight vector of
$B(\infty)$ with weight $0$ is denoted by $u_\infty$.
Similarly $\Ue$ has a crystal base denoted by
$B(-\infty)$, and the unique weight vector of
$B(-\infty)$ with weight $0$ is denoted by $u_{-\infty}$.

Let $\psi$ be the ring automorphism of $\U$ that sends
$\qs$, $e_i$, $f_i$ and $q(h)$ to $\qs$, $f_i$, $e_i$ and $q(-h)$.
It gives a bijection
$B(\infty)\simeq B(-\infty)$ by which
$u_\infty$, $\te_i$, $\tf_i$, $\eps_i$, $\vphi_i$, $\wt$
corresponds to $u_{-\infty}$, $\tf_i$, $\te_i$, $\vphi_i$, $\eps_i$, 
$-\wt$.

Let us denote by $\tU$ the modified quantized
universal enveloping algebra $\oplus_{\lam\in P}\U a_\la$
(see \cite{modified}).
Then $\tU$ has a crystal base $B(\tU)$.
As a crystal,
$B(\tU)$ is regular and isomorphic to
$$\bigsqcup_{\lam\in P}B(\infty)\otimes T_\lam\otimes B(-\infty).$$
Here, $T_\lam$ is the crystal consisting of a single element
$t_\lam$ with $\eps_i(t_\lam)=\vphi_i(t_\lam)=-\infty$ 
and $\wt(t_\lam)=\lam$.

Let $*$ be the anti-involution of $\U$ that sends
$q(h)$ to  $q(-h)$, and $\qs$, $e_i$, $f_i$ to themselves.
The involution $*$ of $\U$ induces
an involution $*$ on $B(\infty)$, $B(-\infty)$, $B(\tU)$.
Then $\te_i^*=*\circ\te_i\circ*$, etc.
give another crystal structure on $B(\infty)$, $B(-\infty)$, $B(\tU)$.
We call it {\em the star crystal structure}.
In the case of $B(\tU)$, these two crystal structures are compatible,
and  $B(\tU)$ may be considered as a
crystal over $\g\oplus \g$.
Hence, for example,
$S^*_w$, the Weyl group action on $B(\tU)$
with respect to the star crystal structure
is a crystal automorphism of $B(\tU)$ with respect
to the original crystal structure.
In particular, the two Weyl group actions 
$S_w$ and $S_{w'}^*$ commute with each other.

The formulas concerning with $B(\tU)$ 
are given in Appendix~\ref{table}.

Note that we have always
\eq\label{eq:epspsi}
&&\text{$\eps_i(b)+\sphi_i(b)=\seps_i(b)+\vphi_i(b)\ge0$ 
for any $b\in B(\infty)$.}
\eneq

\subsection{Schubert decomposition of crystal bases}
For $w\in W$ with a reduced expression
$s_{i_1}\cdots s_{i_\ell}$,
we define the subset $B_w(\infty)$ of 
$B(\infty)$ by
\eq
&&B_w(\infty)=\{\tf_{i_1}^{a_1}\cdots
\tf_{i_\ell}^{a_\ell}u_\infty\set
a_1,\ldots,a_\ell\in \Z_{\ge0}\}.
\eneq
Then $B_w(\infty)$ does not depend on the choice of 
a reduced expression.
We refer the reader to \cite{KL} on the details of
$B_w(\infty)$ and its relationship with the Demazure module.

We have (\cite{KL})
\bnum
\item $B_w(\infty)^*=B_{w^{-1}}(\infty)$.
\item If $w'\le w$, then
$B_{w'(}\infty)\subset B_w(\infty)$.
\item If $s_iw<w$, then $\tf_iB_w(\infty)\subset B_w(\infty)$.
\item $\te_iB_w(\infty)\subset B_w(\infty)\sqcup\{0\}$.
\item
If both $b$ and $\tf_ib$ belong to $B_w(\infty)$,
then all $\tf_i^kb$ ($k\ge0$) belong to $B_w(\infty)$.
\enum
Here $\le$ is the Bruhat order.
Set
\[\bB_w(\infty)=B_w(\infty)\setminus\bigr(
\bigcup_{w'<w}B_{w'}(\infty)\bigr).\]
P. Littelmann (\cite{PL}) showed
$$B(\infty)=\bigsqcup_{w\in W}\bB_{w}(\infty).$$
We have
\eq
&&\bB_w(\infty)^*=\bB_{w^{-1}}(\infty).\\[5pt]
&&\mbox{If $s_iw<w$, then 
\qquad$\te_i^\max \bB_w(\infty)\subset \bB_{s_iw}(\infty)$,}
\label{eq:sch eps}\\
&&\hs{110pt}\tf_i\bB_w(\infty)\subset\bB_w(\infty).\nn\\
&&\mbox{In particular, $\eps_i(b)>0$ for any $b\in\bB_w(\infty)$.}\nn
\eneq
Here, we use the notation $\te_i^\max b=\te_i^{\eps_i(b)}b$.

\subsection{Global bases}\label{subsec:global}
Let $A\subset K$ be the subring of $K$ consisting of rational
functions in $\qs$ without pole at $\qs=0$.
Let $-$ be the automorphism of $K$ sending $\qs$ to $\qs^{-1}$.
Set $K_\Q:=\Q[\qs,\qs^{-1}]$.
Let $V$ be a vector space over $K$,
$L_0$ an $A$-submodule of $V$,
$L_\infty$ an $\ol{A}$- submodule, and
$V_\Q$ a $K_\Q$-submodule.
Set $E:=L_0\cap L_\infty\cap V_\Q$.

\Def[\cite{quant}]
We say that $(L_0,L_\infty,V_\Q)$ is {\em balanced}
if each of $L_0$, $L_\infty$ and $V_\Q$
generates $V$ as a $K$ vector space,
and if the following equivalent conditions are satisfied.
\bnum
\item
$E \to L_0/\qs L_0$ is an isomorphism.
\item
$E \to L_\infty/\qs^{-1}L_\infty$ is an isomorphism.
\item
$(L_0\cap V_\Q)\oplus
(\qs^{-1} L_\infty \cap V_\Q) \to V_\Q$
      is an isomorphism.
\item
$A\otimes_\Q E \to L_0$, $\ol{A}\otimes_\Q E \to L_\infty$,
        $K_\Q\otimes_\Q E \to V_\Q$ and $K \otimes_\Q E \to V$
are isomorphisms.
\enum
\endef

Let $-$ be the ring automorphism of $\U$ sending
$\qs$, $q^h$, $e_i$, $f_i$ to $\qs^{-1}$, $q^{-h}$, $e_i$, $f_i$.

Let $\U_\BQ$ be the $K_\Q$-subalgebra of
$\U$ generated by $e_i^{(n)}$, $f_i^{(n)}$
and $\genfrac{\{}{\}}{0pt}{}{q^h}{n}$ ($h\in P^*$).

Let $M$ be a $\U$-module.
Let $-$ be an involution of $M$ satisfying $(au)^-=\bar a\bar u$
for any $a\in\U$ and $u\in M$.
We call in this paper such an involution a {\em bar involution}.
Let $(L,B)$ be a crystal base of an integrable $\U$-module $M$.

Let $M_\BQ$ be a $\U_\BQ$-submodule of $M$
such that 
\eq\label{eq:zform}
&&
\text{$(M_\BQ){}^-=M_\BQ$, 
and $(u-\ol{u})\in (\qs-1)M_\Q$ for every $u\in M_\Q$.}
\eneq

\Def
If $(L,\ol{L},M_\Q)$ is balanced, we say that $M$ has a global basis.
\endef
In such a case, let $G\colon L/\qs L\isoto E:=L\cap \ol{L}\cap M_\Q$ 
be the inverse of $E\isoto L/\qs L$.
Then $\{G(b);b\in B\}$ forms a basis of $M$.
We call this basis a (lower) {\em global basis}.
The global basis enjoys the following properties
(\cite{quant,global}):
\bnum
\item $\ol{G(b)}=G(b)$ for any $b\in B$.
\item
For any $n\in\Z_{\ge0}$,
$\{G(b);\eps_i(b)\ge n\}$ is a basis
of the $K_\Q$-submodule $\sum_{m\ge n}f_i^{(m)}M_\Q$.
\item
for any $i\in I$ and $b\in B$,
we have
\[f_iG(b)=[1+\eps_i(b)]_iG(\tf_ib)+\sum_{b'}F^{i}_{b,b'}G(b').\]
Here the sum ranges over
$b'\in B$ such that $\eps_i(b')>1+\eps_i(b)$.
The coefficient $F^i_{b,b'}$ belongs to
$\qs q_i^{1-\eps_i(b')}\Q[\qs]$.

 Similarly for $e_iG(b)$.
\enum

\section{Extremal weight modules}
\subsection{Extremal vectors}
Let $M$ be an integrable $\U$-module.
A vector $u\in M$ of weight $\lambda\in P$
is called {\em extremal} (see \cite{AK,modified}),
if we can find
vectors
$\{u_w\}_{w\in W}$
satisfying the following properties:
\eq
&&\text{$u_w=u$ for $w=e$,}\\
&&
\hbox{if $\lan h_i,w\lam\ran\ge 0$, then
$e_iu_w=0$ and $f_i^{(\lan h_i,w\lam\ran)}u_w=u_{s_iw}$,}\\
&&\hbox{if $\lan h_i,w\lam\ran\le 0$, then
$f_iu_w=0$ and $e_i^{(-\lan h_i,w\lam\ran)}u_w=u_{s_iw}$.}
\eneq
Hence if such $\{u_w\}$ exists, then it is unique and
$u_w$ has weight $w\lam$.
We denote $u_w$ by $S_wu$.

Similarly, for a vector $b$ of a regular crystal $B$ with weight $\lam$,
we say that $b$
is an extremal vector
if it satisfies the following similar conditions:
we can find vectors
$\{b_w\}_{w\in W}$ such that
\beq
&&b_w=b\quad\hbox{for $w=e$,}\\
&&\hbox{if $\lan h_i,w\lam\ran\ge 0$ then
$\te_ib_w=0$ and $\tf_i^{\lan h_i,w\lam\ran}b_w=b_{s_iw}$,}\\
&&\hbox{if $\lan h_i,w\lam\ran\le 0$ then
$\tf_iv_w=0$ and $\te_i^{-\lan h_i,w\lam\ran}b_w=b_{s_iw}$.}
\endeq
Then $b_w$ must be $S_wb$.

For $\lam\in P$,
let us denote by $V(\lam)$
the $\U$-module
generated by $u_\lam$
with the defining relation that
$u_\lam$ is an extremal vector of weight $\lam$.
This is in fact infinitely many linear relations on $u_\lam$.
We proved in \cite{modified}
\footnote{In \cite{modified}, it is denoted by $V^\max(\lam)$,
because I thought there would be a natural $\U$-module 
whose crystal base is the connected component of $B(\lam)$.}
that $V(\lam)$ has a global crystal base
$(L(\lam),B(\lam))$.
Moreover the crystal $B(\lam)$ is isomorphic to
the subcrystal of
$B(\infty)\otimes t_\lam\otimes B(-\infty)$
consisting of vectors $b$ such that $b^*$ is an extremal vector
of weight $-\lam$.
We denote by the same letter $u_\lam$ the element of $B(\lam)$
corresponding to $u_\lam\in V(\lam)$.
Then $u_\lam\in B(\lam)$ corresponds to 
$u_\infty\otimes t_\lam\otimes u_{-\infty}$.

Note that, for $b_1\otimes t_\lam\otimes b_2
\in B(\infty)\otimes t_\lam\otimes B(-\infty)$ belonging to $B(\lam)$,
one has
\eq\label{eq:eps*}
&&\parbox{24em}{$\seps_i(b_1)\le\max(\lan h_i,\lam\ran,0)$ and
$\sphi_i(b_2)\le\max(-\lan h_i,\lam\ran,0)$ for any $i\in I$.}
\eneq

For any $w\in W$, $u_\lam\mapsto S_{w^{-1}}u_{w\lam}$
gives an isomorphism of $\U$-modules:
\[V(\lam)\isoto V(w\lam).\]
Similarly, 
letting $S^*_w$ be the Weyl group action on $B(\tU)$
with respect to the  star crystal structure
and regarding $B(\lam)$ as a subcrystal of $B(\tU)$,
$S^*_w\colon B(\tU)\isoto B(\tU)$ induces an isomorphism of crystals
\[S^*_w\colon B(\lam)\isoto B(w\lam).\]

For a dominant weight $\lam$,
$V(\lam)$ is an irreducible highest weight module 
of highest weight $\lam$,
and $V(-\lam)$ is an irreducible lowest weight module
of lowest weight $-\lam$.

\subsection{Dominant weights}

\Def
For a weight $\lam\in P$
and $w\in W$,
we say that
$\lam$ is $w$-dominant $($resp. $w$-regular$)$
if $\lan \beta,\lam\ran\ge 0$
$($resp. $\lan \beta,\lam\ran\neq0$$)$
for any $\beta\in \Delta^\re_+\cap w^{-1}\Delta^\re_-$.
If $\lam$ is $w$-dominant and $w$-regular, we say that 
$\lam$ is regularly $w$-dominant.
\endefinition

If $w=s_{i_\ell}\cdots s_{i_1}$ is a reduced expression, 
then we have
$$\Delta^\re_+\cap w^{-1}\Delta^\re_-
=\{s_{i_1}\cdots s_{i_{k-1}}\alpha_{i_k}\set 1\le k\le\ell\}.$$
Hence $\lam$ is $w$-dominant (resp. $w$-regular)
if and only if
\eq
&&\ba{l}
\lan h_{i_k}, s_{i_{k-1}}\cdots s_{i_1}\lam\ran\ge 0\\
\mbox{(resp. $\lan h_k, s_{i_{k-1}}\cdots s_{i_1}\lam\ran\neq0$).}
\ea
\endeq

Conversely one has the following lemma.
\Lemma
For $i_1,\ldots, i_l\in I$,
and a weight $\lam$,
assume that
\eqn
&&\lan h_{i_k},s_{i_{k-1}}s_{i_{k-1}}\cdots s_{i_1}\lam\ran>0
\quad\mbox{for $k=1,\ldots,l$.}
\endeqn
Then $w=s_{i_{l}}\cdots s_{i_1}$
is a reduced expression.
\enlemma
\proof
By the induction on $l$,
we may assume that
$s_{i_{l-1}}\cdots s_{i_1}$ is a reduced expression.
If $l(w)<l$, then
there exists
$k$ with $1\le k\le l-1$ such that
$s_{i_{l-1}}\cdots s_{i_{k+1}}(h_{i_k})=-h_{i_l}$.
Hence
$$\lan h_{i_l},s_{i_{l-1}}\cdots s_{i_1}\lam\ran
=-\lan h_{i_k},s_{i_{k-1}}s_{i_{k-1}}\cdots s_{i_1}\lam\ran<0,$$
which is a contradiction.
\qed

This lemma implies the following lemma.

\Lemma
Let $w_1$, $w_2\in W$ and
let $\lam$ be an integral weight.
If $\lam$ is regularly $w_2$-dominant and $w_2\lam$ is regularly
$w_1$-dominant,
then $\ell(w_1w_2)=\ell(w_1)+\ell(w_2)$ and
$\lam$ is regularly $w_1w_2$-dominant.
Here $\ell\colon W\to\Z$ is the length function.
\enlemma

\Prop\label{prop:btb}
Let $\lam\in P$ and $b_1\in \bB_{w_1}(\infty)$,
$b_2\in \bB_{w_2}(-\infty)$.
If $b:=b_1\otimes t_\lam\otimes b_2$ belongs to $B(\lam)$, then
one has:
\bnum
\item
$\lam$ is regularly $w_1$-dominant and
$-\lam$ is regularly $w_2$-dominant,
\item
$\ell(w_1w_2^{-1})=\ell(w_1)+\ell(w_2)$,
\item One has
\eqn
S_{w_2}^*(b_1\otimes t_\lam\otimes b_2)
&\in& B_{w_1w_2^{-1}}(\infty)\otimes t_{w_2\lam}\otimes u_{-\infty},\\
S_{w_1}^*(b_1\otimes t_\lam\otimes b_2)
&\in& u_{\infty}\otimes t_{w_1\lam}\otimes B_{w_2w_1^{-1}}(-\infty).
\eneqn
More generally if $w_1=w'w''$ with $\ell(w_1)=\ell(w')+\ell(w'')$,
then
$$S_{w''}^*(b_1\otimes t_\lam\otimes b_2)
\in \bB_{w'}(\infty)\otimes t_{w''\lam}\otimes B_{w_2{w''}^{-1}}(-\infty).$$
\enum
\enprop

\proof
Assume $w_1s_i<w_1$.
Then $c:=\eps^*_i(b_1)>0$ by \eqref{eq:sch eps}.
Hence $\lan h_i,\lam\ran\ge c>0$
by \eqref{eq:eps*}.
We have
$\te_i^*{}^\max b_1\in \bar B_{w_1s_i}(\infty)$.
\eq\label{eq:ind:ext*}
b'=S_i^*(b_1\otimes t_\lam\otimes b_2)
&=&(\te_i^*{}^\max b_1)\otimes
t_{s_i\lam}\otimes (\te_i^*{}^{\lan h_i,\lam\ran-c}b_2).
\endeq
Hence, $\lam$ is regularly $w_1$-dominant by the induction on the length of
$w_1$.
The other statement in (i) is similarly proved.

\vs{5pt}
\noindent
(ii) follows from (i) and the preceding lemma.

\vs{5pt}
In \eqref{eq:ind:ext*},
$\te_i^*{}^{\lan h_i,\lam\ran-c}b_2$ 
belongs to $B_{w_2s_i}(-\infty)$, since (ii) implies $w_2s_i>w_2$.
Repeating this, we obtain (iii).
\qed

\section{Affine quantum algebras}\label{sec:aff}

In the sequel we assume that $\g$ is affine.
\subsection{Affine root systems}
Althogh the materials in this subsection are more or less classical,
we shall review the affine algebras in order to fix
the notations.

Let $\Gg$ be an affine Lie algebra, and let $\Gt$ be its Cartan
subalgebra
(assuming that they are defined over $\Q$).
Let $I$ be the index set of simple roots and
let $\alpha_i\in\Gt^*$ be the simple roots
and $h_i\in\Gt$ the simple coroots
($i\in I$).
We choose a Cartan subalgebra $\Gt$
such that $\{\alpha_i\}_{i\in I}$ and $\{h_i\}_{i\in I}$
are linearly independent and $\dim \Gt=\rank \Gg+1$.
Let us set the root lattice and coroot lattice by
$$Q=\oplus_i\BZ\alpha_i\subset\Gt^*\qbox{and}
Q^\vee=\oplus_i\BZ h_i\subset\Gt.$$
Set $Q_\pm=\pm\sum_i\BZ_{\ge 0}\alpha_i$
and $Q_\pm^\vee=\pm\sum_i\BZ_{\ge 0}h_i$.
Let $\delta\in Q_+$ be a unique element satisfying
$\{\lam\in Q\set \lan h_i,\lam\ran=0\hbox{ for every $i$}\}=\BZ\delta$.
Similarly we define $c\in Q^\vee_+$ by
$\{h\in Q^\vee\set \lan h,\alpha_i\ran=0\hbox{ for every $i$}\}=\BZ c$.
We write
\beq
&&\delta=\sum_ia_i\alpha_i\qbox{and}
c=\sum_ia_i^\vee h_i.
\endeq

We take a $W$-invariant non-degenerate symmetric bilinear form 
$(\cdot,\cdot)$ on $\Gt^*$ normalized by
\beq
&&(\delta,\lam)=\lan c,\lam\ran
\qbox{for any $\lam\in\t^*$.}
\endeq
Then this symmetric form has the signature $(\dim \Gt-1,1)$.
We sometimes identify $\Gt$ and $\Gt^*$ by this symmetric form.
By this identification, $\delta$ and $c$ correspond to each other.

We have
\beq
&&a^\vee_i=\dfrac{(\alpha_i,\alpha_i)}{2}a_i.
\endeq

Note that $(\alpha_i,\alpha_i)/2$ takes the values
$1$, $2$, $3$, $1/2$, $1/3$.
Hence we have for each $i$
\eq
\dfrac{(\alpha_i,\alpha_i)}{2}\in\Z&\hbox{or}&
\dfrac{2}{(\alpha_i,\alpha_i)}\in\Z.
\endeq
If $\Gg$ is untwisted, then 
$2/(\alpha_i,\alpha_i)$ is an integer.

Let us set
$\Gt^*_\cl=\Gt^*/\Q\,\delta$
and let $\cl\colon \Gt^*\to \Gt_\cl^*$ be the canonical projection.
We have
\eqn
&&\t^*_\cl\simeq\bigoplus_{i\in I}(\Q\,h_i)^*.
\eneqn
Set $\Gt^{*0}=\{\lam\in\Gt^*\set \lan c,\lam\ran=0\}$
and $\Gt_\cl^{*0}=\cl(\Gt^{*0})\subset\Gt_\cl^*$.
Then $\clwt$ has a positive-definite symmetric form
induced by the one of $\Gt^*$.

\Lemma For any $a\in\Q$,
$$\cl\colon \{\lam\in\Gt^*\set (\lam,\lam)=a\hbox{ and }(\lam,\delta)\neq0\}
\to\Gt^*_\cl\setminus \clwt$$
is bijective.
\enlemma
\proof
Let $\lam\in\Gt^*$ such that $(\lam,\delta)\not=0$.

Setting $\mu=\lam+x\delta$ for $x\in\Q$,
we have 
$(\mu,\mu)=(\lam+x\delta,\lam+x\delta)
=(\lam,\lam)+2x(\lam,\delta)$.
Hence $\lam+x\delta$ has square length $a$ 
if and only if $x=(a-(\lam,\lam))/2(\lam,\delta)$.
\qed

As a corollary we have
\Prop
$\Gt^*$ endowed with an invariant symmetric form as above, 
simple roots and coroots,
is unique up to a canonical isomorphism.
\enprop

\proof
For example, take
$\rho\in\Gt^*$ such that
$\lan h_i,\rho\ran=1$ for any $i$ and $(\rho,\rho)=0$.
The preceding lemma guarantees its existence and its uniqueness.
The $\alpha_i$'s and $\rho$ form a basis of $\Gt^*$.
\qed

In particular, for any Dynkin diagram isomorphism
$\iota$
(i.e. a bijection $\iota\colon I\to I$ such that
$\lan h_{\iota(i)},\alpha_{\iota(j)}\ran=\lan h_i,\alpha_j\ran$),
there exists a unique isomorphism of $\Gt^*$ that sends
$\alpha_i$ to $\alpha_{\iota(i)}$ and leaves the symmetric form invariant.

Let $\Delta\subset\Gt^*$ be the root system of $\Gg$,
and $\Delta^{\re}$ the set of real
roots: $\Delta^{\re}=\Delta\setminus \Z\,\delta$.
For $\beta\in\t^*$ with $(\beta,\beta)\not=0$,
we set $\beta^\vee=2\beta/(\beta,\beta)$.
Then
$\Delta^\vee:=\{\beta^\vee\set \beta\in \Delta^\re\}\cup
(\Z\,c\setminus\{0\})\subset\Gt$ is
the root system for the dual Lie algebra of $\Gg$.
We set $\Delta^\pm=\Delta\cap Q_\pm$.

Let us denote by $\Delta_\cl$ the image of $\Delta^\real$ by $\cl$.
Then $\Delta_\cl$ is a finite subset of $\Gt_\cl^{*0}$, and
$(\Delta_\cl,\Gt_\cl^{*0})$ is a (not necessarily reduced)
root system.
We call an element of $\Delta_\cl$ a {\em classical} root.

Let $\O(\Gt^*)$ be the orthogonal group of
$\Gt^*$ with respect to the invariant symmetric form.
Let $\O(\Gt^*)_\delta$ be the isotropy subgroup of $\delta$,
i.e. $\O(\Gt^*)_\delta=\{g\in \O(\Gt^*)\set g\delta=\delta\}$.
Then there are canonical group homomorphisms
$$\cl\colon \O(\Gt^*)_\delta\to \GL(\Gt_\cl^*)
\qbox{and}
\cl_0\colon \O(\Gt^*)_\delta\to \O(\Gt_\cl^{*0}).$$

The homomorphism
$\cl\colon \O(\Gt^*)_\delta\to \GL(\Gt_\cl^*)$ is injective.

For $\beta\in \Delta^{\real}$, let $s_\beta$ be the corresponding
reflection $\lam\mapsto \lam-\lan \beta^\vee,\lam\ran\beta$.
Let $W$ be the Weyl group, i.e. the subgroup of $\GL(\Gt^*)$
generated by
the $s_\beta$'s.
Since $W\subset \O(\Gt^*)_\delta$,
there are group homomorphisms
$W\to \GL(\Gt^*_\cl)$ and $W\to \O(\Gt^{*0}_\cl)$.

Let us denote by
$W_\cl$ the image of $W\to \O(\Gt^{*0}_\cl)$.
Then $W_\cl$ is the Weyl group of the root system 
$(\Delta_\cl,\Gt_\cl^{*0})$.

For $\xi\in \Gt^{*0}$, we set
\eq
T(\lam)&=&\lam+(\delta,\lam)\xi-(\xi,\lam)\delta
-\dfrac{(\xi,\xi)}{2}(\delta,\lam)\delta.\label{eq:1}
\endeq
Then $T$ belongs to $\O(\Gt^*)_\delta$, and
$T$ depends only on $\cl(\xi)$.
For $\xio\in \clwt$,
let us define $t(\xio)\in \O(\Gt^*)_\delta$
as the right-hand side of (\ref{eq:1}) with $\xi\in \cl^{-1}(\xio)$.

Then,
\eq
&\ &\mbox{$t\colon \clwt\to\Ker\Big(\cl_0\colon \O(\Gt^*)_\delta\to \GL(\clwt)\Big)$
is a group isomorphism.}
\eneq

We have
\beq
&&g\circ t(\xi)\circ g^{-1}=t\Big(\cl_0(g)(\xi)\Big)
\qbox{for $g\in\O(\Gt^*)_\delta$ and $\xi\in\clwt$.}
\endeq

For $\beta\in \Gt^*$ such that $(\beta,\beta)\not=0$,
let us denote by
$s_\beta$ the reflection
$$s_\beta(\lambda)=\lambda-(\beta^\vee,\lambda)\beta\,.$$
Then we have for $\beta\in \Gt^{*0}$ such that $(\beta,\beta)\not=0$,
\beq
s_{\beta-a\delta}s_\beta=t(a\beta^\vee)\,.
\label{eq:ref2}
\endeq

\medskip
There exists $i_0$
such that 
\beq
\hbox{$W_\cl$ is generated by $\{s_i\set i\neq i_0\}$.}
\label{i_0}
\endeq
If $\Gg$ is not isomorphic to $A^{(2)}_{2n}$,
such an $i_0$ is unique up to a Dynkin diagram automorphism
and $(\alpha_{i_0},\alpha_{i_0})=2$, $a_{i_0}=a_{i_0}^\vee=1$.
In the case of $A_{2n}^{(2)}$,
there are two choices of $i_0$, two extremal nodes, and
$(\alpha_{i_0},\alpha_{i_0})=1$ or $4$,
and accordingly $a_{i_0}=2$ or $1$, $a_{i_0}^\vee=1$ or $2$.

For $\alpha\in\Delta^\re$ or $\alpha\in\Delta_\cl$,
we set
\[c_\alpha=\max(1,\dfrac{(\alpha,\alpha)}{2}),\]
and $c_i=c_{\alpha_i}$.
Then we have, for any $\alpha\in\Delta^\re$
\eq
&&\{n\in\Z\set \alpha+n\delta\in\Delta\}=\Z\,c_\alpha.
\eneq

We set
\beq
&&
\ba{l}
\mbox{$Q_\cl=\cl(Q)$, $Q_\cl^\vee=\cl(Q^\vee)$,
$\tQ=Q_\cl\cap Q_\cl^\vee$.}
\ea\endeq
Here $Q^\vee=\sum_{\alpha\in\Delta^\re}\Z\alpha^\vee$.

We have an exact sequence
\beq
&&\begin{array}{cccccccccc}
1&\longrightarrow&\tQ&\xrightarrow{\ t\ }
&W&\xrightarrow{{\ \cl_0\ }}&W_\cl&\longrightarrow&1\,.
\end{array}
\endeq

For any $\alpha\in \Delta^\re$, let $\tilde\alpha$ be the element
in $\tQ\cap\Q_{>0}\cl(\alpha)$ with the smallest length.
We set
\[\tD=\{\tilde\alpha\set \alpha\in \Delta^\re\}.\]
Then $\tD$ is a reduced root system, and
$\tQ$ is the root lattice of $\tD$.

\Remark
Any affine Lie algebra is either untwisted
or the dual of an untwisted affine algebra
or $A^{(2)}_{2n}$.
\enrem
\begin{itemize}
\item[(i)]
If $\Gg$ is untwisted, then
$\tQ=Q_\cl^\vee\subset Q_\cl$, 
$\tD=\cl(\Delta^{\vee\,\re})$, $\tilde\alpha=\alpha^\vee$.
\item[(ii)]
If $\Gg$ is the dual of an untwisted algebra, then
$\tQ=Q_\cl\subset Q_\cl^\vee$,
$\tD=\cl(\Delta^{\re})$, $\tilde\alpha=\alpha$.
\item[(iii)]
If  $\Gg=A^{(2)}_{2n}$, then
$\tQ=Q_\cl=Q_\cl^\vee$,
$\tD=\cl(\Delta^{\re})=\cl(\Delta^{\vee\,\re})$.
For any $\alpha\in\Delta^\re$, one has
$$\tilde\alpha=
\begin{cases}
\cl(\alpha)&\mbox{ if $(\alpha,\alpha)\not=4$,}\\
\cl(\alpha)/2&\mbox{ if $(\alpha,\alpha)=4$.}
\end{cases}$$
Note that $(\alpha-\delta)/2\in \Delta^\re$ if $(\alpha,\alpha)=4$.
\end{itemize}
If $\g\not=A^{(2)}_{2n}$, then
$\tilde\alpha=c_\alpha\alpha^\vee$.

\Prop\label{pro:tlength}
For $\xi\in\tQ$,
$$l(t(\xi))=\sum_{\beta\in\Delta_\cl}(\beta,\xi)_+/c_\beta
=\dfrac{1}{2}\sum_{\beta\in\Delta_\cl}|(\beta,\xi)|/c_\beta
=\sum_{\beta\in\tD}(\beta^\vee,\xi)_+
.$$
Here $a_+=\max(a,0)$.
\enprop

\proof
For $\beta\in\Delta_\cl$, let us denote by $\beta'$
the unique element of $\Delta^+$ such that
$\cl(\beta')=\beta$ and $\beta'-n\delta
\not\in\Delta^+$ for any $n>0$.
Note that
$(\beta,\xi)\in c_\beta\BZ$.
We have
\eqn
&&t(\xi)^{-1}\Delta^-\bigcap\Delta^+=
\{\gamma\in\Delta^+\set \gamma-(\gamma,\xi)\delta\in\Delta^-\},
\eneqn
and $l(t(\xi))$ is the number of elements in this set.
By setting  $\gamma=\beta'+nc_\beta\delta$, it is isomorphic to
\eqn
&&
\{(\beta,n)\in\Delta_\cl\times\BZ\set 
\,\mbox{$n\ge 0$ and
$\beta'+\Big(nc_\beta-(\beta,\xi)\Big)\delta\in \Delta^-$}\}\\
&&\hs{10pt}
=\{(\beta,n)\in\Delta_\cl\times\BZ\set 
0\le n<(\beta,\xi)/c_\beta\}.
\eneqn
Since $(\beta,\xi)/c_\beta$ is an integer,
we have
\beqn
l(t(\xi))=\sum_{\beta\in\Delta_\cl}
((\beta,\xi)/c_\beta)_+.
\endeqn
The other equalities easily follow.
\qed

\Cor\label{cor:invt}
For $\xi\in\tQ$ and $w\in W_\cl$, 
$$l(t(w\xi))=l(t(\xi)).$$
\encor

We choose a weight lattice
$P\subset\t^*$ satisfying
\eq\label{def:weight}
&&
\quad \left\{
\ba{l}\text{$\alpha_i\in P$ and $h_i\in P^*$ for any $i\in I$.}\\[3pt]
\text{For every $i\in I$, there exists $\Lambda_i\in P$ such that
$\lan h_j,\Lambda_i\ran=\delta_{ji}$.}
\ea\right.
\eneq
We set
\eq
&&\text{$P^0=\{\lam\in P\set\lan c,\lam\ran=0\}$,
$P_\cl=\cl(P)\subset\t^*_\cl$, and $P^0_\cl=\cl(P^0)$.}
\eneq
We have
\[P_\cl=\oplus_{i\in I}(\Z h_i)^*.\]

\Lemma
For $\lam\in P^0$ and $\mu\in \tQ$,
the following two conditions
are equivalent.
\begin{tenumerate}
\item
$\lam$ and $\mu$
are in the same Weyl chamber
$($i.e. for any $\alpha\in \Delta^\re$,
$(\cl(\alpha),\mu)>0$ implies $(\alpha,\lam)\ge0$$)$.
\item
$\lam$ is $t(\mu)$-dominant.
\end{tenumerate}
\enlemma
\proof
For $\alpha\in \Delta^\re$, let us take
$\alpha'\in(\alpha+\Z\delta)\cap\Delta^+$ such that
$\cl(\alpha')=\cl(\alpha)$ and $\alpha'-n\delta\not\in\Delta^+$
for any $n\in\Z_{>0}$.
Then for $\alpha=\alpha'+n\delta\in\Delta^+$,
\eqn
&&\ba{lll}
\alpha\in\Delta^+\cap t(\mu)^{-1}\Delta^{-}
&\Leftrightarrow
t(\mu)\alpha
=\alpha-(\alpha,\mu)\delta\\
&\hs{4em}=\alpha'+(n-(\alpha,\mu))\delta\in\Delta^-
\\[3pt]
&\Leftrightarrow0\le n<(\alpha,\mu).
\ea
\endeqn
(i)$\Rightarrow$(ii)\quad
Now assume $\alpha=\alpha'+n\delta\in \Delta^+\cap t(\mu)^{-1}\Delta^{-}$.
Then $0\le n<(\alpha,\mu)$, and (i) implies $(\alpha,\lam)\ge0$

\vs{5pt}
\noindent
(ii)$\Rightarrow$(i)\quad
Assume $(\alpha,\mu)>0$.
Then taking $n=0$, $\alpha'\in\Delta^+\cap t(\mu)^{-1}\Delta^{-}$,
and hence $(\alpha,\lam)=(\alpha',\lam)\ge0$.
\qed

The following lemma is similarly proved.

\Lemma\label{lem:tdom}
For $\lam\in P^0$ and $\mu\in \tQ$,
the following two conditions
are equivalent.
\bnum
\item
For any $\alpha\in \Delta_\cl$,
$(\alpha,\mu)>0$ implies $(\alpha,\lam)>0$,
\item
$\lam$ is regularly $t(\mu)$-dominant.
\enum
\enlemma


Let us choose $i_0\in I$ as in \eqref{i_0},
and let $W_0$ be the subgroup of $W$ generated by 
$\{s_i\set i\in I\setminus\{i_0\}\}$.
Then $W$ is a semidirect product of $W_0$ and $\tQ$.

\Lemma\label{lem:W0}
Let $\xi\in\tQ$ and $w\in W_0$.
If $\xi$ is regularly $w$-dominant then
\eqn
&&l(t(\xi))=l(t(\xi)w^{-1})+l(w).
\endeqn
\enlemma
\proof
We shall prove the assertion by the induction on $l(w)$.
Write $w=s_iw'$ with $w>w'$ and $i\not=i_0$.
Then
$l(t(\xi))=l(t(\xi)w'{}^{-1})+l(w')$.
Hence it is enough to show
$t(\xi)w'{}^{-1}>t(\xi)w'{}^{-1}s_i$,
or equivalently
$t(\xi)w'{}^{-1}\alpha_i\in\Delta^-$.
We have
\eqn
&&t(\xi)w'{}^{-1}\alpha_i
={w'}^{-1}\alpha_i-(w'\xi,\alpha_i)\delta.
\endeqn
Since $(w'\xi,\alpha_i)>0$, the coefficient of $\alpha_{i_0}$ 
in $t(\xi)w'{}^{-1}\alpha_i$
is negative, and hence $t(\xi)w'{}^{-1}\alpha_i$
is a negative root.
\qed

\subsection{Affinization}\label{subsec:aff}
Let $P$ and $P_\cl$ be as in \eqref{def:weight}.
We denote by $\U$ the quantized universal enveloping algebra
with $P$ as a weight lattice.
We denote by $\Us$
the quantized universal enveloping algebra
with $P_\cl$ as a weight lattice.
Hence $\Us$ is a subalgebra of $\U$
generated by the $e_i$'s, the $f_i$'s
and $q^h$ ($h\in d^{-1}(P_\cl)^*$).
When we talk about an integrable $\U$-module
(resp. $\Us$-module), the weight of its element
belongs to $P$ (resp. $P_\cl$).

Let $M$ be a $\Us$-module
with the weight decomposition $M=\oplus_{\lam\in P_\cl}M_\lam$.
We define a $\U$-module $M_\aff$ with a weight decomposition
$M_\aff=\oplus_{\lam\in P}(M_\aff)_\lam$
by
$$(M_\aff)_\lam=M_{\cl(\lam)}.$$
The action of $e_i$ and $f_i$ are defined in an obvious way, so that
the canonical homomorphism
$\cl\colon M_\aff\to M$ is $\Us$-linear.
We define the $\Us$-linear automorphism $z$ of $M_\aff$
with weight $\delta$
by
$(M_\aff)_\lam\isoto M_{\cl(\lam)}=M_{\cl(\lam+\delta)}\isoto
(M_\aff)_{\lam+\delta}$.

Let us choose $0\in I$ satisfying 
\eq
&&\mbox{
$W_\cl$ is generated by $\{s_i;i\neq 0\}$, and
and $a_{0}=1$.}
\eneq
Recall that $\delta=\sum_ia_i\alpha_i$.
When $\g=A^{(2)}_{2n}$, $0$ is the longest simple root.

Choose a section
$s\colon P_\cl\to P$ of $\cl\colon P\to P_\cl$
such that $s(\cl(\alpha_i))=\alpha_i$ for any $i\in I\setminus\{0\}$.
Then $M$ is embedded into $M_\aff$ by $s$ as a vector space.
We have an isomorphism of $\Us$-modules
\eq\label{emb:aff}
M_\aff\simeq K[z,z^{-1}]\otimes M.
\eneq
Here, $e_i\in\Us$ and $f_i\in\Us$ act on the right hand side by
$z^{\delta_{i0}}\otimes e_i$ and $z^{-\delta_{i0}}\otimes f_i$.

Similarly, for a crystal with weights in $P_\cl$,
we can define its affinization
$B_\aff$ by
\eq
&&B_\aff=\bigsqcup_{\lam\in P}B_{\cl(\lam)}.
\eneq
If an integrable $\Us$-module $M$ has a crystal base 
$(L,B)$, then its affinization
$M_\aff$ has a crystal base $(L_\aff,B_\aff)$.

For $a\in K$, we define the $\Us$-module $M_a$ by
\eq
&&M_a=M_\aff/(z-a)M_\aff.
\eneq

\subsection{Simple crystals}
In \cite{AK}, we defined the notion of simple crystals
and studied their properties.
\Definition\label{unext}
We say that a finite regular crystal $B$ (with weights in $P_\cl^0$) is 
a {\em simple crystal}
if $B$ satisfies
\begin{enumerate}
\item There exists $\lam\in P_\cl^0$ such that
the weight of any extremal vector of $B$ is contained
in $W_\cl\lam$.
\item $\sharp(B_\lam)=1$.
\end{enumerate}
\end{definition}

Simple crystals have the following properties
(loc. cit.).

\Lemma
A simple crystal $B$ is connected.
\enlemma

\Lemma\label{lem:tensimp}
The tensor product of simple crystals is also simple.
\enlemma

\Prop\label{prop:simple-irr}
A finite-dimensional integrable $\Us$-module with a simple crystal base
is irreducible.
\enprop

\section{Affine extremal weight modules}\label{sec:ext}

\subsection{Extremal vectors---affine case}
We prove now one of the main results of this paper.
In the sequel we employ the notations
\[\text{$\te_i^\max b=\te_i^{\eps_i(b)}b$,
$\tf_i^\max b=\tf_i^{\vphi_i(b)}b$,
and similarly for $\te_i^*{}^\max$ and $\tf_i^*{}^\max$.}\]

\Theorem\label{th:fund}
For any $\lam\in P^0$, the weight of any extremal vector
of $B(\lam)$ is contained in $\cl^{-1}\cl(W\lam)$.
\entheorem

\proof
We regard $B(\lam)$ as a subcrystal of 
$B(\infty)\otimes t_\lam\otimes B(-\infty)\subset 
B(\tU)$.

We shall show that $\cl(\wt(b))$ and $-\cl(\wt(b^*))$
are in the same $W_\cl$-orbit
whenever $b$ and $b^*$ are extremal vectors.

For any $b_1\otimes t_\lam\otimes b_2$,
we have 
\beqn
\tf_i^{\max}(b_1\otimes t_\lam\otimes b_2)
&=&b'_1\otimes t_\lam\otimes \tf_i^\max b_2\quad\mbox{for some $b'_1$.}
\endeqn
(For the action of $\tf_i^\max$, etc. on $B(\tU)$, see Appendix~\ref{table}.)
Hence, any extremal vector $b\in B(\lam)$
has the form $b_1\otimes t_\lam\otimes u_{-\infty}$
after applying the $\tf_i^\max$'s.

Hence, we may further assume the following conditions on $b$:
\eq
&&\mbox{$b$ has the form $b_1\otimes t_\lam\otimes u_{-\infty}$,}\\[5pt]
&&\parbox{25em}{for any vector of the form
$b'_1\otimes t_\mu\otimes u_{-\infty}$ in 
$\{S_wS_{w'}^*b\set w,w'\in W\}$,
the length of $\wt(b_1')$ is greater than or equal to
the length of $\wt(b_1)$.}
\eneq
Here, the length of $\sum_im_i\alpha_i$ is by the definition $\sum_i|m_i|$.

Take $i\in I$.
We write $\lam_i=\lan h_i,\lam\ran$
and $\wt_i(b_1)=\lan h_i,\wt(b_1)\ran$ for brevity.

Note that we have $\seps_i(b_1)\le\max(\lam_i,0)$.

We shall show $\wt_i(b_1)\ge0$ for every $i$ in several steps.

\vspace{3pt}
\noindent
(1)\ 
The case $\lam_i\le0$ and $\lam_i+\wt_i(b_1)\le0$.
\hb
Since $b_1\otimes t_\lam\otimes u_{-\infty}$ is a lowest weight vector 
in the $i$-string,
one has 
$\vphi_i(b)=\max(\vphi_i(b_1)+\lam_i,0)=0$,
and hence $\vphi_i(b_1)+\lam_i\le0$.
Similarly, $\seps_i(b)=0$ because $b^*$ is a highest weight vector 
in the $i$-string.
Therefore, one has
\beqn
S_i^*S_i(b_1\otimes t_\lam\otimes u_{-\infty})
&=&\tf_i^*{}^{-\lam_i}(\te_i{}^\max b_1\otimes t_{\lam}
\otimes \te_i^{-\vphi_i(b_1)-\lam_i}u_{-\infty})\\
&=&(\tf_i^*{}^{\vphi_i(b_1)}\te_i{}^\max b_1)\otimes t_{s_i\lam}
\otimes u_{-\infty}.
\endeqn
The last equality follows from
$S_i^*S_i(b)=(\tf_i^*{}^{k}\te_i{}^\max b_1)\otimes t_{s_i\lam}
\otimes u_{-\infty}$ for some $k$.

Hence, the minimality of $b_1$ gives
$$0\le\varphi_i(b_1)-\eps_i(b_1)=\wt_i(b_1).$$

\smallskip
\noindent
(2)\ 
The case $\lam_i>0$ and $\lam_i+\wt_i(b_1)\le0$.
\hb
We shall show that this case cannot occur.
In this case, as in (i), 
\beqn
\varphi_i(b_1)+\lam_i\le0.
\endeqn
On the other hand,
$\varphi^*_i(b_1\otimes t_\lam\otimes u_{-\infty})
=\max(\eps^*_i(b_1)-\lam_i,0)=0$
implies 
\beqn
\eps^*_i(b_1)\le\lam_i.
\endeqn
Hence we obtain (the first inequality by \eqref{eq:epspsi})
\beqn
&&0\le\eps^*_i(b_1)+\varphi_i(b_1)
=(\eps^*_i(b_1)-\lam_i)+(\varphi_i(b_1)+\lam_i)
\le0,
\endeqn
which implies
$\eps^*_i(b_1)=\lam_i$
and
$\varphi_i(b_1)=-\lam_i$.
Then
we have
\beqn
&&\te_i^*{}^\max(b_1\otimes t_\lam\otimes u_{-\infty})
=(\te^*_i{}^\max b_1)\otimes t_{s_i\lam}\otimes u_{-\infty}.
\endeqn
Hence, the minimality of $\wt(b_1)$ implies $\eps^*_i(b_1)=0$, and
this contradicts $\eps^*_i(b_1)=\lam_i>0$.

\smallskip
\noindent
(3)\ 
The case $\lam_i\ge0$ and $\lam_i+\wt_i(b_1)\ge0$.
\hb
In this case, one has $\eps_i(b)=\sphi_i(b)=0$, and hence
$\vphi_i(b)=\lam_i+\wt_i(b_1)$,
which implies $\varphi_i(b)-(\lam_i-\seps_i(b_1))=\sphi_i(b_1)\ge0$.
Hence we have
\eqn
S_iS_i^*(b_1\otimes t_\lam\otimes u_{-\infty})
&=&\tf_i^{\vphi_i(b)}(\te_i^*{}^\max b_1\otimes t_{s_i\lam}\otimes
\te_i^{\lam_i-\seps_i(b_1)}u_{-\infty})\\
&=&(\tf_i^{\sphi_i(b)}\te^*_i{}^\max b_1)\otimes t_{s_i\lam}\otimes
u_{-\infty}.
\eneqn
Hence we have
$\varphi_i^*(b_1)\ge \eps_i^*(b_1)$,
or equivalently $\wt_i(b_1)\ge0$.

\smallskip
\noindent
(4)\ 
The case $\lam_i\le0$ and $\lam_i+\wt_i(b_1)\ge0$.\hb
We have immediately
$\wt_i(b_1)\ge0$.

\medskip
In all the cases we have $\wt_i(b_1)\ge0$.
Since $\wt(b_1)$ is of level $0$,
one has
$0=\lan c,\wt(b_1)\ran=\sum_ia^\vee_i\wt_i(b_1)$,
which implies that
$\wt_i(b_1)=0$ for every $i$, or equivalently
$\cl(\wt(b_1))=0$.
\qed

\Cor
For any $\lam\in P$,
the weight of any vector in $B(\lam)$
is contained in the convex hull of
$W\lam$.
\encor
\proof
In the positive level case (i.e. $\lan c,\lam\ran>0$), 
$\lam$ being conjugate to a dominant
weight and $B(\lam)$ is isomorphic to
the crystal base of an irreducible highest weight module.
In this case, the assertion is well-known.
Similarly for negative level case.

Assume that the level of $\lam$ is zero. 
Note that all vector in $B(\lam)$ can be reached at an extremal vector
after applying $\te_i^\max$ and $\tf_i^\max$
by \cite{modified}.
Hence the assertion follows
from the preceding theorem.
Note that
$\cl^{-1}\cl(W\lam)$ is contained in the convex hull of
$W\lam$ provided that $\cl(\lam)\not=0$.
\qed

The following theorem is an immediate consequence of the preceding corollary.

\Theorem
Let $M$ be an integrable $\Us$-module
and $u$ a vector in $M$ of weight $\lam\in P_\cl$.
Then the following conditions are equivalent.
\bnum
\item
$u$ is an extremal vector.
\item
The weights of $\Us u$ are contained in the convex hull of $W_\cl\lam$.
\item
$\Us_\beta u=0$ for any $\beta\in\Delta_\cl$ such that
$(\beta,\lam)\ge0$.
\enum
\entheorem

In particular, for any $\lam\in P$,
$V(\lam)$ is isomorphic to the $\U$-module
generated by a weight vector $u$ of weight $\lam$ with
(iii) in the above corollary and the following integrability
condition as defining relations:
$$\text{$f_i^{1+\lan h_i,\lam\ran}u=0$ if $\lan h_i,\lam\ran\ge0$ and
$e_i^{1-\lan h_i,\lam\ran}u=0$ if $\lan h_i,\lam\ran\le0$.}$$
%
%
%

\subsection{Fundamental representations}\label{ss:fd}
Let us take $\dz\in I$ such that
\eq\label{eq:0hat}
&&\mbox{
$W_\cl$ is generated by $\{s_i;i\neq \dz\}$, and
and $a_{\dz}^\vee=1$.}
\eneq
Recall that $c=\sum_ia_i^\vee h_i$.
When $\g=A^{(2)}_{2n}$, $\dz$ is the shortest simple root.
We set $I_{\dz}=I\setminus\{\dz\}$.
For $i\in I_{\dz}$, we set 
$$\varpi_i=\Lambda_i-a_i^\vee\Lambda_{\dz}\in P^0.$$
Hence we have
$P^0_\cl=\oplus_{i\in I_\dz}\Z\cl(\varpi_i)$.
We say that $\lam\in P$ is a 
{\em \basic\ weight}
if $\cl(\lam)$ is $W_\cl$-conjugate to some $\cl(\varpi_i)$
($i\in I_{\dz}$).
Note that this notion does not depend on the choice
of $\dz$.
\Prop
Assume that $\lam=\sum_{i\in J}\varpi_i$
for some subset $J$ of $I_\dz$.
Then one has:
\bnum
\item
any extremal vector of $B(\lam)$ is in the $W$-orbit of
$u_\lam$,
\item
$B(\lam)$ is connected.
\enum
\enprop
\proof
(ii) follows from (i)
because any vector is connected with extremal vector.

Let us prove (i).
We use arguments similar to the proof of Theorem~\ref{th:fund}.
Let us take an extremal vector $b\in B(\lam)$.
Among the vectors
in $S_wS_{w'}^*b$ with the form 
$b_1\otimes t_\mu\otimes u_{-\infty}$,
we take one such that $\wt(b_1)$ has the smallest length.
Then the proof in Theorem~\ref{th:fund} shows that
$\cl(\wt(b_1))=0$.
Hence, one has
\beqn
S_iS^*_i(b_1\otimes t_\mu\otimes u_{-\infty})
&=&
\begin{cases}
\tf_i{}^{\seps_i(b_1)}\te^*_i{}^\max(b_1)
\otimes t_{s_i\mu}\otimes u_{-\infty}
&\mbox{if $\mu_i\ge0$,}\\
\tf^*_i{}^{\eps_i(b_1)}\te_i{}^\max(b_1)
\otimes t_{s_i\mu}\otimes u_{-\infty}
&\mbox{if $\mu_i\le0$.}
\end{cases}
\endeqn
In the both cases, the length of $b_1$ remains unchanged
after applying $S_iS^*_i$.
Therefore, applying $S_{w'{}^{-1}}S^*_{w'{}^{-1}}$,
we can assume $w'=1$ and $\mu=\lam$.

For $i\in I\setminus J$, we have $\lam_i\le0$, which implies
$\eps^*_i(b_1)=0$.
If $i\in J$, then $\lam_i=1$ and hence $\eps^*_i(b_1)$ ($\le\lam_i$)
must be $0$ or $1$. On the other hand, we have
\beqn
&&S^*_i(b_1\otimes t_\lam\otimes u_{-\infty})=
\te^*_i{}^\max b_1\otimes t_\lam\otimes
\te_i^{\lam_i-\eps^*_i(b_1)}u_{-\infty}.
\endeqn
If $\eps^*_i(b_1)=1$, then
this contradicts the minimality of $\wt(b_1)$.
Hence $\eps^*_i(b_1)=0$ for every $i\in J$.

Thus we have $\eps^*_i(b_1)=0$ for every $i\in I$ and hence
$b_1=u_\infty$.
Thus we obtain $u_\lam=S_wb$.
\qed

The following theorem is a particular case of the preceding proposition.
\Theorem\label{th:bas_mul1}
If $\lam\in P$ is a \basic\ weight,
then
any extremal vector of $B(\lam)$ is in the $W$-orbit of
$u_\lam$.
\entheorem
We shall now study further properties of $B(\lam)$ for a \basic\
weight $\lam$.

\Lemma\label{lem:iso}
Let $\lam$ be a \basic\ weight.
Then $\{w\in W;w\lam=\lam\}$ is generated by
$\{s_\beta\set\beta\in \Delta^\re_+\,,\,(\beta,\lam)=0\}$.
\enlemma
\proof
We may assume $\lam=\Lambda_j-a_j^\vee\Lambda_\dz$
for some $j\in I_\dz$.
Since the similar statement holds for $(W_\cl,\t_\cl^{*0}$),
it is enough to show that
$t(\xi)$ is contained in the subgroup $G$ generated by
$\{s_\beta\set\beta\in \Delta^\re_+\,,\,(\beta,\lam)=0\}$,
provided that $\xi\in\tQ$ and $(\xi,\lam)=0$.
We have
$s_{a\delta-\beta}s_{\beta}=t(a\beta^\vee)$
by (\ref{eq:ref2}).
In particular, one has
$t(c_\beta\beta^\vee)\in G$
whenever $\beta\in\Delta^\re$ satisfies $(\beta,\lam)=0$.

\vs{4pt}
\noindent
(1) The case where $\g\not=A^{(2)}_{2n}$\quad
It is enough to show that
$\{\xi\in\tQ;(\xi,\lam)=0\}$
is generated by
$\{c_\beta\beta^\vee;
\beta\in\Delta_\cl,\,(\beta,\lam)=0\}$.
In this case, $\tQ$ has a basis
$\{c_i\alpha_i^\vee\set i\in I_\dz\}$.
Hence
$\{\xi\in\tQ\set(\xi,\lam)=0\}$
is generated by $\{c_i\alpha_i^\vee\set i\in I_\dz\setminus\{j\}\}$.

\vs{4pt}
\noindent
(2) The case where $\g=A^{(2)}_{2n}$\quad
In this case, $\tQ=Q=\oplus_{i\in I_\dz}\Z\tilde\alpha_i$.
Hence
$\{\xi\in\tQ\set(\xi,\lam)=0\}$
has a basis $\{\tilde\alpha_i\set i\in I_\dz\setminus\{j\}\}$.
Hence, the result follows from
\eqn
&&t(\tilde\alpha_i)=
\begin{cases}
s_{\delta-\alpha_i}s_{\alpha_i}
&\mbox{if $(\alpha_i,\alpha_i)=2$,}\\
s_{(\delta-\alpha_i)/2}s_{\alpha_i}
&\mbox{if $(\alpha_i,\alpha_i)=4$.}
\end{cases}
\eneqn
Note that $(\delta-\alpha_i)/2$ is a real root in the last case.
\qed

\Lemma\label{lem:refb}
For any $\beta\in\Delta^\re$ and any $\lam\in P$ 
such that $s_\beta\lam=\lam$,
we have
$S_{s_\beta}(u_\infty\otimes t_\lam\otimes u_{-\infty})=
u_\infty\otimes t_\lam\otimes u_{-\infty}$.
\enlemma
\proof
Set $a_\lam=u_\infty\otimes t_\lam\otimes u_{-\infty}$.
We assume $\beta\in\Delta^\re_+$.
We shall prove the assertion by the induction on the length of
$\beta$.
If $\beta$ is a simple root, it is obvious.
Otherwise, we can write $\beta=s_i\gamma$ 
for a positive real root $\gamma$ 
whose length is less than that of $\beta$.
We have
$S_{s_\beta}=S_iS_{s_\gamma}S_i$.
Set $\mu=s_i\lam$. Then
$s_\gamma\mu=\mu$ and hence we have
$S_{s_\gamma}a_\mu=a_\mu$ by the induction hypothesis.
Since $S_iS_i^*a_\lam=a_\mu$
or equivalently $S_ia_\lam=S_i^*a_\mu$,
we have
\[
S_{s_\beta}a_\lam=S_iS_{s_\gamma}S_ia_\lam
=S_iS_{s_\gamma}S^*_ia_\mu=
S_iS^*_ia_\mu=a_\lam.
\]
\qed

Lemma~\ref{lem:iso} and Lemma~\ref{lem:refb}
imply the following proposition.
\Prop\label{prop:ind_weyl}
Let $\lam$ be a \basic\ weight.
\bnum
\item
If $w\in W$ satisfies $w\lam=\lam$, then $S_wu_\lam=u_\lam$ and
$S^*_wu_\lam=u_\lam$.
\item
For $\mu\in W\lam$, the isomorphism
$S^*_w\colon B(\lam)\isoto B(\mu)$
does not depend on $w\in W$ such that $\mu=w\lam$.
\enum
Here we regard $B(\lam)$ and $B(\mu)$ as subcrystals of $B(\tU)$.
\enprop

\Remark
For a general $\lam\in P^0$, it is not true
that the extremal weights of $B(\lam)$ belong to $W\lam$.
For example in $\lam=2(\Lambda_1-\Lambda_0)$ in the $A^{(1)}_1$-case
$f_0f_1u_\lam$ is an extremal vector with weight $\lam-\delta$.
\enrem
\Remark
It is not true in general $w\lambda=\lambda$
implies $S_wu_\lambda=u_\lambda$.
For example in the case of
$\Gg=A^{(1)}_2$, and $\lambda=\Lambda_1+\Lambda_2-2\Lambda_0$,
set $w_1=t(\alpha_1)=s_1s_0s_2s_1$ and $w_2=t(\alpha_2)=s_1s_0s_1s_2$.
Then $w_1\lambda=w_2\lambda=\lambda-\delta$,
but $S_{w_1}u_\lambda\not=S_{w_2}u_\lambda$.
\enrem
\Conj
For any $\lam\in P$,
$S_w u_\lam=u_\lam$ if and only if
$w\in W$ is in the subgroup generated by
$\{s_\beta; \mbox{$\beta$ is a real root such that $(\beta,\lam)=0$}\}$.
\enconj
%

Theorem~\ref{th:bas_mul1}
and Proposition~\ref{prop:ind_weyl}
immediately imply the following result.

\Prop
Assume that $\lam$ is a \basic\ weight.
\bnum
\item
$B(\lam)_\lam=\{u_\lam\}$.
\item
$B(\lam)$ is connected.
\enum
\enprop
\proof
Let $b\in B(\lam)_\lam$.
Then Theorem~\ref{th:bas_mul1} implies $b=S_wu_\lam$ for
some $w\in W$ with $w\lam=\lam$, and
Proposition~\ref{prop:ind_weyl} implies
$S_{w}u_\lam=u_\lam$.
\qed

In order to show the finite multiplicity theorem
for $B(\varpi_i)$, we shall need the following result.
\Lemma\label{lem:circ}
Assume $\lam=\cl(\Lam_{i_1}-a^\vee_{i_1}\Lam_0)$ 
for some $i_1\in I_{\dz}$
and $\mu\in W_\cl\lam$.
If $w\in W$ satisfies $l(w)\ge \sharp W_\cl$ and
$\mu$ is regularly $w$-dominant,
then there exist $w'$, $w''\in W$ such that
$w=w'w''$, $l(w)=l(w')+l(w'')$ and
$\lam=w''\mu$.
\enlemma
\proof
Let $w=s_{i_1}\cdots s_{i_l}$ be a reduced expression of $w$.
Since $l\ge\sharp W_\cl$, there exists
$0\le j<k\le l$ such that
$\cl(s_{i_1}\cdots s_{i_j})=\cl(s_{i_1}\cdots s_{i_k})$.
Hence $s_{i_{j+1}}\cdots s_{i_k}=t(\xi)$
for some $\xi\in\tQ\setminus\{0\}$.
Replacing $\mu$ with $s_{i_{k+1}}\cdots s_{i_\ell}\mu$,
we reduce the lemma to the following sublemma.
\qed
\Sublemma 
If $\xi\in\tQ\setminus\{0\}$ and $\mu\in W_\cl\lam$ is regularly $t(\xi)$-dominant,
then there exists $w_1\in W$ such that
$\lam=w_1\mu$ and $l(t(\xi))=l(t(\xi)w_1^{-1})+l(w_1)$.
\ensublemma
\proof
Let us take $w\in W_{\dz}:=\lan s_i;i\in I_\dz\ran$ such that
$\mu=w\lam$ and $\lam$ is regularly $w$-dominant.
By Lemma~\ref{lem:tdom}, for $\beta\in \Delta_\cl$,
$(\beta,\xi)>0$ implies $(\beta, \mu)>0$.
Hence $(\beta, w^{-1}\xi)>0$ implies $(\beta, \lam)>0$.
In particular, $(\beta, \lam)=0$ (resp. $(\beta, \lam)>0$)
implies $(\beta, w^{-1}\xi)=0$ (resp. $(\beta, w^{-1}\xi)\ge0$).
For $i\in I_\dz\setminus\{i_1\}$,  $(\alpha_i, w^{-1}\xi)=0$
because $(\alpha_i,\lam)=0$.
Moreover $(\alpha_{i_1}, w^{-1}\xi)\ge 0$
because $(\alpha_{i_1},\lam)>0$.
Hence we have
$w^{-1}\xi=c\lam$ for $c>0$.
Hence $w^{-1}\xi$ is regularly $w$-dominant.
Corollary~\ref{cor:invt} and Lemma~\ref{lem:W0} imply that
\eqn
l(t(\xi))&=&l(t(w^{-1}\xi))=l(t(w^{-1}\xi)w^{-1})+l(w)
=l(w)+l(w^{-1}t(\xi)).
\endeqn
Then the sublemma follows by setting
$w_1=w^{-1}t(\xi)$.
\qed

\Prop
Let $\lam\in P$ be a \basic\ weight.
Then for every $\xi\in P$, $B(\lam)_\xi$ is a finite set.
\enprop
\proof
For $w\in W$ and $\mu\in W\lam$, we define a subset $A_w(\mu)$ of
$B(\tU)$
by
\[A_w(\mu)=\{b\otimes t_\mu\otimes u_{-\infty}\in B(\mu);
b\in\bB_{w}(\infty)\},\]
and then set
$A_w=\bigsqcup_{\mu\in W\lam}A_w(\mu)$.
Note that $A_w(\mu)$ is a finite set.
One has
$$B(\lam)\subset \bigcup_{w,w_1\in W}S^*_{w_1}(A_{w}(w_1^{-1}\lam)).$$
We shall first show
\eq\label{eq:fundres}
B(\lam)\subset 
\bigcup_{\substack{w_1\in W,\\w\in W
\ with\ \ell(w)\le N}}S^*_{w_1}(A_{w}).
\eneq
Here $N:=\sharp W_\cl$.

For $b:=b_1\otimes t_\mu\otimes u_{-\infty}$ in $A_{w}$,
we shall show
$$b\in \bigcup_{\substack{w_1\in W,\\w'\in W
\ with\ \ell(w')\le N}}S^*_{w_1}(A_{w'})$$
by the induction on $\ell(w)$.

Proposition~\ref{prop:btb} implies that $\mu$ is regularly $w$-dominant.
We may assume $\ell(w)>N$.
By Lemma~\ref{lem:circ}, there exists
$w_1=w'w''$ such that $l(w)=l(w')+l(w'')$, $w'\neq 1$ 
and $\lam':=w''\mu$ satisfies $\cl(\lam')=\cl(\lam)$.

By Proposition~\ref{prop:btb}, one has
\eqn
&&S_{w''}^*(b_1\otimes t_\mu\otimes u_{-\infty})
=b'_1\otimes t_{\lam'}\otimes b'_2
\endeqn
with
$b_1'\in\bB_{w'}(\infty)$
and $b_2'\in B_{{w''}^{-1}}(-\infty)$.
Take $i\in I$ such that $w's_i<w'$. Then $\lam'_i>0$ implies $i=i_1$.
Hence $c:=\eps_i^*(b_1')\le \lam'_i=1$.
One has
\eq
S_i^*(b_1'\otimes t_{\lam'}\otimes b'_2)
=(\te_i^*{}^\max b'_1)\otimes t_{s_i\lam'}
\otimes \te_i^*{}^{\lam'_i-c}b'_2.
\endeq
If $c=1$, then $\lam'_i-c=0$.
Take $x\in W$ such that
$b_2'\in\bB_{x}(-\infty)$.
Then $x\le {w''}^{-1}$, since $b_2'\in B_{{w''}^{-1}}(-\infty)$.
Since $\te_i^*{}^\max b'_1\in \bB_{w's_i}(\infty)$,
Proposition~\ref{prop:btb} implies
\[S^*_x((\te_i^*{}^\max b'_1)\otimes t_{s_i\lam}
\otimes b'_2)
\in B_{w's_ix^{-1}}(\infty)\otimes t_{xs_i\lam}\otimes u_{-\infty}.\]
Since $\ell(w's_ix^{-1})<\ell(w)$, the induction proceeds.

Next assume $c=0$. Then $\lam_j\le0$ for $j\in I\setminus\{i_1\}$
implies $\seps_j(b'_1)=0$ for every $j\in I$.
Hence $b_1'=u_\infty$. This contradicts
$w'\not=1$ and $b_1'\in\bB_{w'}(\infty)$.
Thus we have proved \eqref{eq:fundres}.

\vs{5pt}
For $\mu\in W\lam$, set
$$C(\mu)=
\bigcup_{w\in W\ with\ \ell(w)\le N}A_{w}(\mu).$$
Taking $w\in W$ such that $\mu=w\lam$, we set
$$\tilde C(\mu):=S^*_{w^{-1}}C(\mu)\subset B(\lam),$$
By Proposition~\ref{prop:ind_weyl},
$\tilde C(\mu)$ does not depend on the choice of $w$.
We have
\bnum
\item
$\tilde C(\mu)$ is a finite set,
\item
there is a finite subset $F$ of $Q$ independent of $\mu$
such that $\Wt(\tilde C(\mu))\subset \mu+F$.
\enum
Hence, for any $\xi\in P$,
\[B(\lam)_\xi
\subset
\bigcup_{\mu\in W\lam}\tilde C(\mu)_\xi
=\bigcup_{\mu\in W\lam\cap(\xi-F)}\tilde C(\mu)_\xi\]
is a finite set.
\qed

We have thus obtained the following properties of $V(\lam)$.
\Prop\label{prop:fundrep}
Let $\lam\in P^0$ be a \basic\ weight.
\bnum
\item
$\Wt(V(\lam))$ is contained in the intersection of 
$\lam+Q$ and the convex hull of $W\lam$.
\item
$\dim V(\lam)_\mu=1$ for any $\mu\in W\lam$.
\item
$\dim V(\lam)_\mu<\infty$ for any $\mu\in P$.
\item
$\Wt(V(\lam))\cap(\lam+\Z\,\delta)\subset W\lam$.
\item
$V(\lam)$ is an irreducible $\U$-module.
\item
Any non-zero integrable $\U$-module generated by an extremal weight
vector of weight $\lam$ is isomorphic to $V(\lam)$.
\enum
Moreover $V(\lam)$ has a global base.
\enprop
For any $\mu\in W\lam$, let us denote by
$u_\mu$ the unique global basis in $V(\lam)_\mu$.
Since $u_\mu$ is an extremal vector with weight $\mu$,
we have the $\U$-linear homomorphism
$V(\mu)\to V(\lam)$
that sends $u_\mu\in V(\mu)$ to $u_\mu\in V(\lam)$.
This homomorphism is in fact an isomorphism.

Set $\lam=\varpi_i$. One has
\eq
&&
\{n\in \Z;\varpi_i+n\delta\in W\varpi_i\}=\Z d_i,
\endeq
where $d_i=(\varpi_i,\tilde\alpha_i)$.
Note that $d_i=\max(1,(\alpha_i,\alpha_i)/2)\in\Z$
except the case $d_i=1$ when $\g=A^{(2)}_{2n}$ and $\alpha_i$ is the
longest root.
Hence one has
$$\bigoplus_{\mu\in\cl^{-1}\cl(\varpi_i)}V(\varpi_i)_\mu
=\bigoplus_{n\in\Z}V(\varpi_i)_{\lam+nd_i}.$$
We have a $\U$-linear isomorphism
$V(\varpi_i+d_i\delta)\isoto V(\varpi_i)$.
Since there is a $\Us$-linear isomorphism
$V(\varpi_i)\isoto V(\varpi_i+d_i\delta)$
that sends $u_{\varpi_i}$ to $u_{\varpi_i+d_i\delta}$,
we obtain a $\Us$-linear automorphism $z_i$
of $V(\varpi_i)$ of weight $d_i\delta$,
which sends $u_{\varpi_i}$ to $u_{\varpi_i+d_i\delta}$.

Let us define the $\Us$-module $W(\varpi_i)$ by
\eq
&&W(\varpi_i)=V(\varpi_i)/(z_i-1)V(\varpi_i).
\eneq

The following result is now obvious.
\begin{theorem}\label{th:fundrep}
\bnum
\item
$W(\varpi_i)$ is a finite-dimensional irreducible integrable
$\Us$-module.
\item
$W(\varpi_i)$ has a global basis with a simple crystal.
\item
For any $\mu\in \Wt(V(\varpi_i))$,
$$W(\varpi_i)_{\cl(\mu)}\simeq V(\varpi_i)_\mu.$$
\item
$\dim W(\varpi_i)_{\cl(\varpi_i)}=1$.
\item
The weight of any extremal vector of $ W(\varpi_i)$ 
belongs to $W\cl(\varpi_i)$.
\item $\Wt(W(\varpi_i))$ is the intersection of
$\cl(\varpi_i)+Q_\cl$ and the convex hull of
$W\cl(\varpi_i)$.
\item
$K[z_i^{1/d_i}]\otimes_{K[z_i]}V(\varpi_i)\simeq W(\varpi_i)_\aff$.
Here the action of $z_i^{1/d_i}$ on the left hand side corresponds to
the action of $z$ on the right hand side
defined in \S~\ref{subsec:aff}.
\item
$V(\varpi_i)$ is isomorphic to the submodule
$K[z^{d_i},z^{-d_i}]\otimes W(\varpi_i)$
of $W(\varpi_i)_\aff$ as a $\U$-module.
Here we identify $W(\varpi_i)_\aff$ 
with $K[z,z^{-1}]\otimes W(\varpi_i)$
as in \eqref{emb:aff}.
\item
Any irreducible finite-dimensional integrable $\Us$-module
with $\cl(\varpi_i)$ as an extremal weight
is isomorphic to $W(\varpi_i)_a$ for some $a\in K\setminus\{0\}$.
\enum
\end{theorem}
\proof
The irreducibility of $W(\varpi_i)$ follows for example by
Proposition~\ref{prop:simple-irr},
and the other assertions are now obvious.
\qed

We call $W(\varpi_i)$ a  {\em fundamental representation} (of level $0$).

\section{Existence of Global bases}\label{sec:exis}
\subsection{Regularized modified operators}

For $n\in\Z$ and $i\in I$,
let us define the operator $\TF$
\eq
&&\TF=\sum_{k\ge0, -n}f_i^{(n+k)}e_i^{(k)}a_k(t_i).
\eneq
Here 
\[a_k(t_i)=(-1)^kq_i^{k(1-n)}t_i^k\prod_{\nu=0}^{k-1}(1-q_i^{n+2\nu}).\]
Then it acts on any integrable $\U$-module $M$. 
Moreover it acts also on
any $\U_\BQ$-submodule $M_\BQ$.
In this sense, $\TF$ has no pole except $q=0$, $\infty$.
Let $(L,B)$ be a crystal base of $M$.
Then we have the following result,
which says that $\TF$ has no pole at $q=0$
and coincides with $\tf^n_i$ at $q=0$.

\Prop
We have $\TF L\subset L$,
and the action of $\TF$ on $L/\qs L$
coincides with $\tf_i^n$.
\enprop
\proof
In order to prove this, it is sufficient to prove the following
statement.
For any weight vector $u\in M$ with $e_iu=0$ and $m\in\Z_{\ge0}$, we have
\[
\TF f_i^{(m)}u=cf_i^{(m+n)}u\]
for some $c\in K:=\Q(\qs)$ regular at $\qs=0$ and $c(0)=1$.
Set $t_iu=q_i^lu$.
Then we can assume
\[l\ge n+m.\]
We have
\[a_k(t_i)f_i^{(m)}u=a_k(q_i^{l-2m})f_i^{(m)}u.\]
Hence
\eqn
f_i^{(m)}u
&=&\sum_{k\ge0}a_k(q_i^{l-2m})f_i^{(n+k)}e_i^{(k)}f_i^{(m)}u\\
&=&\sum_{k=0}^ma_k(q_i^{l-2m})f_i^{(n+k)}\qbin{l-m+k}{k}f_i^{(m-k)}u\\
&=&\sum_{k=0}^ma_k(q_i^{l-2m})\qbin{n+m}{m-k}\qbin{l-m+k}{k}f_i^{(m+n)}u.
\eneqn
Here, $$\qbin{n}{m}=\dfrac{[n]_i!}{[m]_i![n-m]_i!}$$
is the $q$-binomial coefficient.
Hence it is enough to show that
\eqn
&&A:=\sum_{k=0}^ma_k(q_i^{l-2m})\qbin{n+m}{m-k}\qbin{l-m+k}{k}
\in 1+q_i\Z[q_i].
\eneqn
This follows immediately from the following formula,
whose proof due to Anne Schilling is
given in Appendix~\ref{app:qhg}.
\eq\label{eq:anne}
&&A=
\sum_{k=0}^{m}q_i^{k(2l-2m-n+2)}
\prod_{j=1}^k\dfrac{1-q_i^{n+2(j-1)}}{1-q_i^{2j}}
\prod_{j=1}^{m-k}\dfrac{1-q_i^{n+2j}}{1-q_i^{2j}}.
\eneq
\qed

\subsection{Existence theorem}
We shall use the notations and terminologies in \S~\ref{subsec:global}.
Let $M$ be an integrable $\U$-module,
$-$ a bar involution of $M$, and
$(L,B)$ a crystal base of $M$.
Let $M_\BQ$ be a $\U_\BQ$-submodule of $M$
such that $(M_\BQ){}^-=M_\BQ$.
Set $E:=L\cap \overline L\cap M_\BQ$.

\Theorem\label{th:ex}
Let $S$ be a subset of $P$.
We assume the following conditions:
\begin{tenumerate}
\item
$\{(\xi,\xi);\xi\in \Wt(M)\}$ is bounded from above.
\item
$u-\bar u\in (\qs-1)M_\BQ$ for any $u\in M_\BQ$.
\item
$M_\BQ$ generates $M$ as a vector space over $K$.
\item
For any $\xi\in P\setminus S$, $(L_\xi, \ol L_\xi, (M_\BQ)_\xi)$
is balanced.
\item
Any extremal weight $($i.e. the weight of an extremal vector$)$
of $B$ is in $P\setminus S$.
\item
$\qs L\cap \ol L\cap M_\BQ=0$.
\end{tenumerate}
Then we have
\begin{anumerate}
\item
$(L, \ol L, M_\BQ)$
is balanced.
\item
For any $n$, we have
\[f_i^nM=\bigoplus\limits_{\eps_i(b)\ge n}\BQ(\qs)G(b)\quad\mbox{and}\quad
e_i^nM=\bigoplus\limits_{\vp_i(b)\ge n}\BQ(\qs)G(b).\]
\item
$M_\BQ=\sum_{\xi\in P\setminus S}\U_\BQ(M_\BQ)_\xi$ and
$M=\sum_{\xi\in P\setminus S}\U M_\xi$.
\end{anumerate}
\end{theorem}

\vs{1em}
The rest of this section is devoted to the proof of this theorem.

\Lemma\label{lem:1.1}
The action of $-$ on $E$ is the identity.
\enlemma
\proof
For $u\in E$, we have
$(u-\bar u)/(1-\qs^{-1})\in \qs L\cap \ol L\cap M_\BQ=0$.
\qed

By (vi), the homomorphism
$E\to L/\qs L$ is injective.
Let us denote by $B'$ the intersection of $B$ 
and the image of this homomorphism.
To see (a), it is enough to show that $B=B'$.
For $b\in B'$, let us denote by $G(b)$
the element $E$ such that $b\equiv G(b)\ \mod\, \qs L$.
Note that $G(b)\ol{\phantom {a}}=G(b)$ by Lemma~\ref{lem:1.1}.
We shall prove 
the following statements by the descending induction on $(\xi,\xi)$:
{\allowdisplaybreaks
\eq
&&
\parbox{25em}
{$B_\xi=B'_\xi$,
or equivalently, $(L_\xi,\, \ol L_\xi,\, (M_\BQ)_\xi)$ is
balanced,}
\label{BB'}\\[5pt]
&&
\parbox{25em}
{$G(b)-f_i^{(\eps_i(b))}G(\te^\max_ib)
\in \sum_{\eps_i(b')>\eps_i(b)}\BQ[\qs,\qs^{-1}]G(b')$
for any $b\in B_\xi$,}\label{BB'1}
\\[6pt]
&&
\begin{array}{l}
\sum\limits_{b\in B_\xi,\,\eps_i(b)\ge n}\BQ[\qs,\qs^{-1}]G(b)
=\sum_{m\ge n}f_i^{(m)}(M_\BQ)_{\xi+m\alpha_i}\qquad\\[-3pt]
\hfill\text{for any $n\ge\max(0,-\lan h_i,\xi\ran)$.}
\end{array}
\label{BB'2}
\endeq
}
as well as the similar statements replacing $f_i$ with $e_i$.


If $(\xi,\xi)$ is big enough, those statements are trivially satisfied
by (i).
Now assuming \eqref{BB'}--\eqref{BB'2} for
$\xi$ such that $(\xi,\xi)>a$, let us prove
them for $\xi$ with $(\xi,\xi)=a$.

\Lemma\label{lem:1.2}
Let $i\in I$.
Set $k=\max(0,-\lan h_i,\xi\ran)$.
\banum
\item
If $\te^\max_ib\in B'$, then $b\in B'$ and
\[G(b)-f_i^{(\eps_i(b))}G(\te^\max_ib)
\in \sum_{\substack{b'\in B'_\xi\\\eps_i(b')>\eps_i(b)}}
\BQ[\qs,\qs^{-1}]G(b').\]
In particular, any $b\in B_\xi$ with $\eps_i(b)>k$ 
is contained in $B'$.
\item
$\sum\limits_{\substack{b\in B'_\xi\\\eps_i(b)\ge n}}\BQ[\qs,\qs^{-1}]G(b)
=\sum_{m\ge n}f_i^{(m)}(M_\BQ)_{\xi+m\alpha_i}$
for any $n>k$.
\eanum
The similar statements hold after exchanging $e_i$ and $f_i$.
\enlemma
\proof
Let us prove the lemma by the descending induction on $n$
(in the case (a), $n$ means $\eps_i(b)$, and hence $n\ge k$).
If $n$ is big enough, they are true by the hypothesis (i) in Theorem~
\ref{th:ex}.
Let us prove (a).
Set $b_1=\te^\max_ib$.
Then
$u=\TF G(b_1)$ satisfies
$b\equiv u\ \mod\, \qs L$ and
\[u-f_i^{(n)}G(b_1)\in 
\sum_{m>0}\Z[\qs,\qs^{-1}]f_i^{(m+n)}e_i^{(m)}G(b_1)
\subset\sum_{m>n}f_i^{(m)}(M_\BQ)_{\xi+m\alpha_i}.\]
The induction hypothesis (b) implies that the last space is contained in 
\[\sum_{\substack{b'\in B'_\xi\\ \eps_i(b')>n}}\Q[\qs,\qs^{-1}]G(b').\]
Hence we can write
$u-f_i^{(n)}G(b_1)=\sum_{b'} c_{b'}G(b')$ 
where $b'$ ranges over $b'\in B'$ with $\eps_i(b')>n$
and $c_{b'}\in\BQ[\qs,\qs^{-1}]$.
Hence we can write $c_{b'}-\ol{c_{b'}}=c'_{b'}-\ol{c'_{b'}}$
with $c'_{b'}\in \qs\BQ[\qs]$.
Then $v:=u-\sum_{b'} c'_{b'}G(b')
=f_i^{(n)}G(b_1)+\sum_{b'} (c_{b'}-c'_{b'})G(b')$
satisfies $\ol v=v$ and hence it belongs to $E$.
Moreover one has $b\equiv v\ \mod\, \qs L$.
Hence $b$ belongs to $B'$, and $G(b)=v$.

To complete the proof of (a),
it is enough to remark
$\te_i^\max b\in B'$ when $\eps_i(b)>k$, because
$(\wt(\te^\max_i(b)),\,\wt(\te^\max_i(b)))>(\wt((b),\wt(b))$.

Let us prove (b).
The left hand side is contained in the right hand side 
by (a) and the induction hypothesis on $n$.
Let us show the opposite inclusion.
Set $\eta=\xi+n\alpha_i$ with $n>k$. Then we have $(\eta,\eta)>(\xi,\xi)$,
and \eqref{BB'2} holds for $\eta$.
Hence we have
$$(M_\Q)_\eta\subset\sum_{\eps_i(b)=0,\,b\in B'_\eta} \BQ[\qs,\qs^{-1}]G(b)
+\sum_{m>0}f_i^{(m)}(M_\Q)_{\eta+m\alpha_i},$$
which implies 
\eqn
f_i^{(n)}(M_\BQ)_{\eta}&\subset&
\sum_{\eps_i(b)=0,\,b\in B'_\eta} \BQ[\qs,\qs^{-1}]f_i^{(n)}G(b)
+\sum_{m>n}f_i^{(m)}M_\BQ\\
&\subset&\sum_{\substack{\eps_i(b)=0\\b\in B'_\eta}} \BQ[\qs,\qs^{-1}]f_i^{(n)}G(b)
+\sum_{\substack{\eps_i(b)>n\\b\in B'_\xi}} \BQ[\qs,\qs^{-1}]G(b).
\eneqn
The desired inclusion follows from (a).

\qed

\Lemma\label{lem:aux:mt}
$B_\xi\subset B'$.
\enlemma
\proof
Let $b\in B_\xi$.
By the hypothesis (v), 
there exists $X_l\cdots X_1b$ whose weight is outside $S$,
where $X_\nu$ is $\te_i^\max$ or $\tf_i^\max$.
Hence by the induction on $l$ we
may assume that $\te_i^\max b$ or $\tf_i^\max b$ is contained in $B'$.
Then the preceding lemma implies $b\in B'$.
\qed

The properties \eqref{BB'} and \eqref{BB'1} are now obvious,
and \eqref{BB'2} easily follows from
Lemma \ref{lem:1.2} and Lemma\ref{lem:aux:mt}.

Thus the induction proceeds,
and we complete the proof of (a), (b) in Theorem~\ref{th:ex}.

Finally let us prove (c).
Set $M'=\sum_{\xi\in P\setminus S}\U M_\xi$ and
$M'_\BQ=\sum_{\xi\in P\setminus S}\U_\BQ(M_\BQ)_\xi$.
Set $L'=L\cap M'$.
Then $L'$ is invariant by $\te_i$ and $\tf_i$.
By the hypothesis (v), any vector in $B$ 
is connected with a vector whose weight is outside
$S$.
Hence $B$ is contained in $L'/\qs L'\subset L/\qs L$.
This shows that $(L', B)$ is a crystal base of $M'$,
and $L'/\qs L'=L/\qs L$.
Thus we can apply Theorem~\ref{th:ex} to $M'$.
Hence we obtain
$L'\cap\ol{L'}\cap M'_\Q=L\cap\ol{L}\cap M_\Q$,
and
$M'_\Q=K_\Q\otimes(L'\cap\ol{L'}\cap M'_\Q)=M_\Q$.
This completes the proof of Theorem~\ref{th:ex}.

\section{Universal $R$-matrix}\label{sec:Rmatrix}
In this section, we shall review the universal $R$-matrix
introduced by Drinfeld
and the universal bar involution introduced by Lusztig.

Although we mainly use the following coproduct 
$\Delta$ in this article 
\eq
\begin{array}{rl}
\Delta(q^h)&=q^h\otimes q^h\\
\Delta(e_i)&=e_i\otimes t_i^{-1}+1\otimes e_i\\
\Delta(f_i)&=f_i\otimes 1+t_i\otimes f_i,
\end{array}
\eneq
we shall introduce another coproduct $\Db=(-\otimes-)\circ\Delta\circ-$
\eq
\begin{array}{rl}
\Db(q^h)&=q^h\otimes q^h\\
\Db(e_i)&=e_i\otimes t_i+1\otimes e_i\\
\Db(f_i)&=f_i\otimes 1+t_i^{-1}\otimes f_i.
\end{array}
\eneq

Let $M_\nu$ ($\nu=1,2$) be
a $\U$-module with weight decomposition.
Let us denote by $M_1\otimes M_2$ the tensor product of $M_1$ and
$M_2$
with the $\U$-module structure induced by $\Delta$,
and $M_1\bt M_2$ the $\U$-module induced by $\Db$.

Then there is an isomorphism
\[q^{-(\,\cdot\,,\,\cdot\,)}\colon M_1\bt M_2\to M_2\otimes M_1\]
given by
\[q^{-(\,\cdot\,,\,\cdot\,)}\,(x\bt y)=q^{-(\wt(x),\wt(y))}\,y\otimes x.\]

Let us define the ring $\Up\ct\Um$ by
\eq
&&\Up\ct\Um=\bigoplus_{\xi\in
  Q}\prod_{\xi=\lam+\mu}\bigl(\Up_\lam\otimes\Um_\mu\bigr).
\eneq
The counits $\Up\to K$ and $\Um\to K$ induces
$\eps\colon \Up\ct\Um\to K$.
 Modifying Drinfeld's construction  (\cite{D0})
of a universal R-matrix, Lusztig has shown that
there exists a unique intertwiner
$\Xi\in\CT$ satisfying the following properties:
\[\text{$\Xi\circ\Delta(a)=\Db(a)\circ\Xi$ for any $a\in\U$,}\]
normalized by $\eps(\Xi)=1$.
Then it satisfies
\eq\label{eq:inv:xi}
\ol{\Xi}\circ \Xi=\Xi\circ\ol{\Xi}=1.
\eneq

We introduce the completion of the tensor products as follows.
We set
\eqn
&&F_{(\lam,\mu)}(M_1\ct M_2)=\prod
_{\gamma\in Q_+}
(M_1)_{\lam+\gamma}\otimes (M_2)_{\mu-\gamma}\\
&&F_{>(\lam,\mu)}(M_1\ct M_2)=\prod
_{\gamma\in Q_+\setminus\{0\}}
(M_1)_{\lam+\gamma}\otimes (M_2)_{\mu-\gamma},
\eneqn
and then
\[
M_1\ct M_2=\sum_{\lam,\mu\in P}
F_{(\lam,\mu)}(M_1\ct M_2)\subset \prod_{\lam,\mu\in
  P}(M_1)_\lam\otimes (M_2)_\mu.\]

Sometimes we use another completion $M_1\rct M_2$ in the opposite direction:
\eqn
&&F_{(\lam,\mu)}(M_1\rct M_2)=\prod
_{\gamma\in Q_+}
(M_1)_{\lam-\gamma}\otimes (M_2)_{\mu+\gamma}
\eneqn
and then
\[
M_1\rct M_2=\sum_{\lam,\mu\in P}
F_{(\lam,\mu)}(M_1\rct M_2)\subset \prod_{\lam,\mu\in
  P}(M_1)_\lam\otimes (M_2)_\mu.\]

They have a structure of a $\U$-module
by $\Delta$ and containing $M_1\otimes M_2$ as a $\U$-submodule.

We denote by
$M_1\cbt M_2$ the same vector space $M_1\ct M_2$
with the action of $\U$ induced by
$\Db$. Then $M_1\cbt M_2$ contains $M_1\bt M_2$ as a $\U$-submodule.

We have an isomorphism
\[q^{-(\,\cdot\,,\,\cdot\,)}\colon M_1\cbt M_2\isoto M_2\rct M_1.\]
The operator $\Xi$ induces an isomorphism
\[M_1\ct M_2\isoto M_1\cbt M_2.\]
Then $\Xi$ sends $F_{(\lam,\mu)}(M_1\ct M_2)$ to
$F_{(\lam,\mu)}(M_1\cbt M_2)$, and
\eq\label{eq:xi}
&&
\begin{array}{l}
\mbox{The homomorphism induced by $\Xi$}\\[3pt]
\quad M_1{}_\lam\otimes M_2{}_\mu\simeq
F_{(\lam,\mu)}(M_1\ct M_2)/F_{>(\lam,\mu)}(M_1\ct M_2)\\[3pt]
\quad\qquad\longrightarrow
F_{(\lam,\mu)}(M_1\cbt M_2)/F_{>(\lam,\mu)}(M_1\cbt M_2)
\simeq M_1{}_\lam\otimes M_2{}_\mu\\[3pt]
\mbox{is equal to the identity.}
\end{array}\eneq

The intertwiner
$\Ru\colon M_1\ct M_2\to M_2\rct M_1$, called
the {\em universal $R$-matrix},
is given by
by
\eq\label{eq:univR}
&&\Ru\colon M_1\ct M_2\xrightarrow{\ \Xi\ }
M_1\cbt M_2\xrightarrow{q^{-(\,\cdot\,,\,\cdot\,)}}
M_2\rct M_1.
\eneq
It is an isomorphism.

Assume that $M_1$ and $M_2$ have a bar involution.
Then \eqref{eq:inv:xi} implies that
\[\cju\colon M_1\ct M_2\xrightarrow{\ \Xi\ }M_1\cbt M_2
\xrightarrow{-\otimes -}M_1\ct M_2\]
is a bar involution on $M_1\ct M_2$ as observed by 
G. Lusztig (\cite{GL}).
We call it the {\em universal bar involution}.

\section{\Good\ modules}\label{sec:good}
Let us take a finite-dimensional integrable $\Us$-module $M$.
We consider the following conditions on $M$:
\eq
&&\mbox{$M$ has a bar involution,}\\
&&\mbox{$M$ has a crystal base $(L(M),B(M))$,}\\
&&\mbox{$M$ has a global base,}\\
&&\mbox{$B(M)$ is a simple crystal.}
\eneq
In this paper, we say that a $\Us$-module $M$ is a {\em \good}
$\Us$-module
if $M$ satisfies the above conditions.
The level zero fundamental representations $W(\vpi_i)$
is a \good\ $\Us$-module.
A \good\ $\Us$-module is always irreducible 
(Proposition~\ref{prop:simple-irr}).

Let $M_1$ and $M_2$ be \good\ $\Us$-modules.
Then we have
\eqn
(M_1)_\aff\ct (M_2)_\aff&=&K[[z_1/z_2]]\bigotimes_{K[z_1/z_2]}
\Bigl((M_1)_\aff\otimes (M_2)_\aff\Bigr),\\
(M_2)_\aff\rct (M_1)_\aff&=&K[[z_1/z_2]]\bigotimes_{K[z_1/z_2]}
\Bigl((M_2)_\aff\otimes (M_1)_\aff\Bigr).
\eneqn
Here $z_\nu$ is the $\Us$-linear automorphism
of weight $\delta$ on $(M_\nu)_\aff$
introduced in \S~\ref{subsec:aff}.
\Lemma\label{lem:genirr}
$K(z_1/z_2)\otimes_{K[z_1/z_2]}\bigl((M_1)_\aff\otimes (M_2)_\aff\bigr)$ 
is an irreducible module
over $K(z_1/z_2)\otimes_{K[z_1/z_2]}\U[z_1^{\pm1},z_2^{\pm1}]$.
\enlemma
\proof
Since $M_1\otimes M_2$ has a simple crystal base
by Lemma~\ref{lem:tensimp},
it is irreducible by Proposition~\ref{prop:simple-irr}.
Then the lemma follows from the fact that the specialization of 
$(M_1)_\aff\otimes (M_2)_\aff$ at the special point $z_1/z_2=1$
is irreducible.
\qed

By the result of the previous section, we have
the bar involution
$$\cju\colon (M_1)_\aff\ct (M_2)_\aff\to(M_1)_\aff\ct (M_2)_\aff.$$
It commutes with $z_1$ and $z_2$.
Let $u_\nu$ be the extremal vector with dominant weight
$\lam_\nu$ of $M_\nu$ ($\nu=1,2$), and set $u=u_1\otimes u_2$.
Then we have
$\bigl((M_1)_\aff\ct (M_2)_\aff\bigr)_{\lam_1+\lam_2}
=K((z_1/z_2))u$.
Hence, by \eqref{eq:xi}, we have
\eq\label{eq:phi}
&&\text{$\cju(u)=\ol{\vphi(z_1/z_2)}u\ $ 
or equivalently $\ \Xi(u)=\vphi(z_1/z_2)u$}
\eneq
for some $\vphi(z_1/z_2)\in K[[z_1/z_2]]$ with $\vphi(0)=1$.
We define 
\[\cjn\colon (M_1)_\aff\ct (M_2)_\aff\to(M_1)_\aff\ct (M_2)_\aff\]
by $\cjn=\cju\circ\vphi(z_1/z_2)^{-1}$.
Then it satisfies
$$\cjn(u)=u.$$
\Lemma
$\cjn$ is a unique endomorphism of 
$(M_1)_\aff\ct(M_2)_\aff$ satisfying
$\cjn(u_1\otimes u_2)=u_1\otimes u_2$ and
$\cjn(av)=\ol{a}\cjn(v)$ for any $a\in \U((z_1/z_2))[z_2^{\pm1}]$,
$v\in (M_1)_\aff\ct(M_2)_\aff$.
\enlemma
\proof
It is enough to show that
a $\U[z_1^{\pm1},z_2^{\pm1}]$-linear homomorphism
$$f\colon (M_1)_\aff\otimes (M_2)_\aff\to(M_1)_\aff\ct
(M_2)_\aff$$
vanishes if $f(u_1\otimes u_2)=0$.
By Lemma~\ref{lem:genirr},
$K(z_1/z_2)\otimes_{K[z_1/z_2]}(M_1)_\aff\otimes (M_2)_\aff$ 
is an irreducible module over
$K(z_1/z_2)[z_2^{\pm1}]\otimes\U$.
Hence the assertion follows.
\qed
\noindent

Hence, $\cjn$ defines a bar involution on $(M_1)_\aff\ct (M_2)_\aff$,
which we call the {\em normalized bar involution}.
In particular we have
\[\vphi(z)\ol{\vphi(z)}=1.\]
In the sequel, we use the normalized
bar involution to define a global basis.

The universal $R$-matrix:
\[\Ru\colon (M_1)_\aff\ct (M_2)_\aff\to(M_2)_\aff\rct (M_1)_\aff\]
sends $u_1\otimes u_2$ to
$q^{-(\lam_1,\lam_2)}\vphi(z_1/z_2)u_2\otimes u_1$ 
with the same function $\vphi$ given in \eqref{eq:phi}.
Hence setting
$\Rn=q^{(\lam_1,\lam_2)}\vphi(z_1/z_2)^{-1}\Ru$, we have an intertwiner
\[\Rn\colon (M_1)_\aff\ct (M_2)_\aff\to(M_2)_\aff\rct (M_1)_\aff\]
that sends $u_1\otimes u_2$ to
$u_2\otimes u_1$. We call $\Rn$ the {\em normalized $R$-matrix}.
Both $R$-matrices commute with $z_1$ and $z_2$. 

By \eqref{eq:univR}
and \eqref{eq:xi},
we have, for any $v_\nu\in (M_\nu)_\aff$, 
\begin{multline}\label{eq:normR}
\Rn(v_1\otimes v_2)
\equiv q^{\lan \lam_1,\lam_2\ran-\lan\wt(v_1),\wt(v_2)\ran}v_2\otimes v_1\\
\mod \prod_{\xi\in Q_+\setminus\{0\}}
((M_2)_\aff)_{\wt(v_2)-\xi}
\otimes ((M_1)_\aff)_{\wt(v_1)+\xi}.
\end{multline}
We have also
\eqn
\Rn\colon (M_1)_\aff\otimes (M_2)_\aff\
&\to& K(z_1/z_2)\otimes_{K[z_1/z_2]}\bigl((M_2)_\aff\otimes(M_1)_\aff\bigr)\\
&&\hspace*{10pt}\hookrightarrow (M_2)_\aff\rct(M_1)_\aff.
\eneqn

We shall generalize these observations to the case of
tensor products of several modules.
Let $M_\nu$ ($\nu=1,\ldots,m$) be a \good\ $\Us$-modules
with a crystal base $(L_\nu,B_\nu)$.
Let $(M_\nu)_\aff$ be its affinization.
Then $(M_\nu)_\aff$ has a crystal base 
$\bigl((L_\nu)_\aff,(B_\nu)_\aff\bigr)$.
Let $\lam_\nu\in P$ be a dominant extremal weight of $(M_\nu)_\aff$,
and $u_\nu$ the extremal global basis with weight $\lam_\nu$.
We denote the canonical automorphism $(M_\nu)_\aff$
of weight $\delta$ by $z_\nu$.

Then 
\[M:=\otimes_{\nu=1}^m(M_\nu)_\aff=(M_1)_\aff\otimes\cdots\otimes(M_m)_\aff\]
has a structure of $K[z_1^{\pm1},\ldots, z_m^{\pm1}]$-module.
Set
\eq
M&=&(M_1)_\aff\otimes\cdots\otimes (M_m)_\aff,\\
M_\Q&=&(M_1{}_\Q)_\aff\otimes\cdots\otimes (M_m{}_\Q)_\aff,
\eneq
and let $\bigl(L(M), B(M)\bigr)$ be the tensor product of
the crystal bases of the $(M_\nu)_\aff$'s.
We set
\[
\Mh=
K[[z_1/z_2,\ldots,z_{m-1}/z_m]]\mathop\otimes\limits_{K[z_1/z_2,\ldots,z_{m-1}/z_m]}
\otimes_{\nu=1}^m(M_\nu)_\aff.\]
We set also

\eqn
L(\Mh)&=&
A[[z_1/z_2,\ldots,z_{m-1}/z_m]]
\mathop\otimes\limits_{A[z_1/z_2,\ldots,z_{m-1}/z_m]}
\otimes_{\nu=1}^m(L_\nu)_\aff,\\
\Mh_\Q&=&
\Q[[z_1/z_2,\ldots,z_{m-1}/z_m]]
\mathop\otimes\limits_{\Q[z_1/z_2,\ldots,z_{m-1}/z_m]}
\otimes_{\nu=1}^m((M_\nu)_\aff{}_\Q).
\eneqn

Similarly to the case of the tensor product of two modules,
we can define the universal bar involution of $\Mh$
by $$\cju=(-\otimes\cdots\otimes-)\circ\prod_{1\le i<j\le m}\Xi_{ij},$$
where
$\Xi_{ij}$ is the operator $\Xi$ acting on the $i$-th and $j$-th
components of the tensor product.
Normalizing $\cju$, we obtain the normalized bar involution $\cjn$ on
$\Mh$.
It satisfies, by setting $u=u_1\otimes\cdots\otimes u_m$,
\[\cjn(u)=u.
\]
Moreover it satisfies for $v_\nu\in (M_\nu)_\aff$
\mult\label{eq:cjn:tri}
\cjn(v_1\otimes \cdots\otimes v_m)\equiv
\ol{v_1}\otimes\cdots\otimes\ol{v_m}\\[5pt]
\mod \prod_{\xi_1,\ldots,\xi_m}
((M_1)_\aff)_{\wt(v_1)+\xi_1}\otimes\cdots\otimes
((M_m)_\aff)_{\wt(v_m)+\xi_m}.
\end{multline}
Here the product ranges over
$\xi_1,\dots,\xi_m\in Q$
with
$\sum_{\nu=1}^m\xi_\nu=0$
and $\sum_{\nu=1}^\mu\xi_{\nu}\in Q^+\setminus\{0\}$
($\mu=1,\ldots,m-1$).

Since $\cjn$ is expressed by a triangular matrix,
the well-known argument of triangular matrices implies the following result.
\Lemma\label{lem:bal}
\bnum
\item
$\qs L(\Mh)\bigcap\bigl(\cjn L(\Mh)\bigr)\bigcap \Mh_\Q=0$.
\item For any $b=b_1\otimes\cdots b_m\in B(M)$,
there exists a unique
$G(b)\in L(\widehat M)$ such that
$\cjn(G(b))=G(b)$ and $b\equiv G(b)\ \mod\ \qs  L(\widehat M)$.
\item
Moreover $G(b)$ has the form
\[\hs{40pt}G(b)=
G(b_1)\otimes\dots\otimes G(b_m)
+\sum c_{b_1',\cdots,b_m'}G(b'_1)\otimes\dots\otimes G(b'_m).
\]
Here the infinite sum ranges over 
$b_1'\otimes\cdots\otimes b'_m\in B(M)$
such that
$\sum\limits_{\nu=1}^m\wt(b'_\nu)=\sum\limits_{\nu=1}^m\wt(b_\nu)$
and
$\sum\limits_{\nu=1}^\mu(\wt(b'_\nu)-\wt(b_\nu))\in Q_+\setminus\{0\}$
$(\mu=1,\dots,m-1)$.
Moreover $c_{b_1',\cdots,b_m'}\in \qs\Q[\qs]$.
\enum
\enlemma

Later we shall see that this infinite sum is in fact a finite sum.

Set 
\[N=\U[z_1^{{\pm1}},\dots,z_m^{{\pm1}}]u.\]
Then $N$ is a submodule of $\Mh$ stable by the bar involution $\cjn$.
Set $\lam=\sum_{\nu=1}^m\lam_\nu$.
Then we have
\eqn
N_{\lam+\Z\delta}
&:=&\oplus_{n\in\Z}N_{\lam+n\delta}
=\bigl(\otimes_{\nu=1}^m(M_\nu)_\aff\bigr)_{\lam+\Z\delta}\\
&=&\otimes_{\nu=1}^m\bigl((M_\nu)_\aff\bigr)_{\lam_\nu+\Z\delta}
=K[z_1^{{\pm1}},\dots,z_m^{{\pm1}}](u_1\otimes\cdots\otimes u_m).
\eneqn
Hence one has
\eq
&&\text{$N_\mu=M_\mu$ for any $\mu\in W\lam+\Z\delta$.}
\eneq
Define 
\eqn
L(N)&=&L(M)\cap N,\\
N_\Q&=&M_\Q\cap N,\\
B(N)&=&B(M).
\eneqn
Then $L(N)/\qs L(N)\subset L(M)/\qs L(M)$.

\Lemma\label{lem:Ncr}
$B(N)$ is a basis of $L(N)/\qs L(N)$,
and $(L(N),B(N))$ is a crystal base of $N$.
\enlemma
\pf
Since $B(N)$ is a basis of $L(M)/\qs L(M)$,
it is enough to show that $B(N)$ 
is contained in $L(N)/\qs L(N)$.
Since every vector in $B(N)$ is connected with an extremal vector with
weight in
$\lam+\Z\delta$,
and extremal vectors with such a weight is
$u$ up to the action of $z_1^{\pm1},\ldots,z_m^{\pm1}$,
we obtain the desired result.
\qed

Setting $S=\Wt(M)\setminus(W\lam+\Z\delta)$, we can apply
Theorem~\ref{th:ex} to $N$. The hypotheses in the theorem
are satisfied by  Lemma~\ref{lem:bal} and Lemma~\ref{lem:Ncr}, 
and we obtain the following theorem.

\Theorem\label{th:glex}
\begin{tenumerate}
\item
$(L(N),\cjn L(N), N_\Q)$ is balanced.
Hence $N$ has a global base.
\item
$N_\Q=\U_\Q[z_1,\ldots z_m]u$.
\end{tenumerate}
\entheorem
Furthermore, Lemma \ref{lem:bal} implies the following proposition.
\Prop\label{prop:exp:glob}
For any $b_\nu\in B((M_\nu)_\aff)$ $(\nu=1,\dots m)$,
we have
\[G(b_1\otimes\dots \otimes b_m)=
G(b_1)\otimes\dots\otimes G(b_m)
+\sum c_{b_1',\cdots,b_m'}G(b'_1)\otimes\dots\otimes G(b'_m).\]
Here the sum ranges over 
$(b_1',\dots,b'_m)\in \prod_{\nu=1}^mB((M_\nu)_\aff)$
such that
$\sum_{\nu=1}^m\wt(b'_\nu)=\sum_{\nu=1}^m\wt(b_\nu)$
and
$\sum_{\nu=1}^\mu(\wt(b'_\nu)-\wt(b_\nu))\in Q_+\setminus\{0\}$
$(\mu=1,\dots,m-1)$.
Moreover $c_{b_1',\cdots,b_m'}\in \qs\Q[\qs]$,
and $c_{b_1',\cdots,b_m'}$ vanishes except for finitely many
$(b_1',\cdots,b_m')$.
\enprop

By specializing at $z_\nu=1$, we obtain the following proposition.
\Prop
The tensor product of \good\ $\Us$-modules 
is also a \good\ $\Us$-module.
\enprop

\section{Main theorem}\label{sec:main}
The following theorem is conjectured in \cite{AK}
in the special case when all the $M_\nu$ are fundamental representations.
Note that, as seen by the proof,
the theorem holds
even if we consider $\U$ as an algebra over the algebraically closed field
$\widehat K:=\sum\limits_{n>0} \C((q^{1/n}))$, and
$a_\nu$ as elements of $\widehat K$, and replace $A$ with
the subring $\widehat A:=\sum\limits_{n>0} \C[[q^{1/n}]]$ of $\widehat K$.

\Theorem
\bi
\item[(i)]
Let $M_\nu$ {\rm ($\nu=1,\ldots,m$)} be \good\ $\Us$-modules.
Let $a_\nu\in K$.
Assume that $a_{\nu}/a_{\nu+1}\in A$ for $\nu=1,\dots,m-1$.
Then $(M_1)_{a_1}\otimes (M_2)_{a_2}\otimes \cdots\otimes
(M_m)_{a_m}$
is generated by $u_1\otimes\dots\otimes u_m$.
\item[(ii)]
Assume that $(M_\nu)^*$ {\rm ($\nu=1,\ldots,m$)}
is a \good\ $\Us$-module,
and $a_{\nu+1}/a_{\nu}\in A$ for $\nu=1,\dots,m-1$.
Then any non-zero 
submodule of $(M_1)_{a_1}\otimes (M_2)_{a_2}\otimes \cdots\otimes
(M_m)_{a_m}$
contains $u_1\otimes\dots\otimes u_m$.
\ei
\end{theorem}

\proof
Since (ii) is the dual statement of (i),
it is enough to prove (i).

Let us embed the crystal $B_\nu$ of $M_\nu$
into $(B_\nu)_\aff$ as in $\eqref{emb:aff}$.
Let $\psi\colon (M_1)_{\aff}\otimes \cdots\otimes(M_m)_{\aff}
\to(M_1)_{a_1}\otimes \cdots\otimes(M_m)_{a_m}$
be the canonical projection.
Then $\psi(G(b_1)\otimes\cdots\otimes G(b_m))$ ($b_\nu\in B_\nu$)
forms a basis of $(M_1)_{a_1}\otimes (M_2)_{a_2}\otimes \cdots\otimes
(M_m)_{a_m}$.
Since
$\psi(G(b_1\otimes\cdots\otimes b_m))$ 
are in $\Us (u_1\otimes\dots\otimes u_m)$, it is enough to show that
they also generate
$(M_1)_{a_1}\otimes (M_2)_{a_2}\otimes \cdots\otimes
(M_m)_{a_m}$ as a vector space.

By Proposition~\ref{prop:exp:glob},
we can write
\begin{multline*}
G(b_1\otimes\dots \otimes b_m)=
G(b_1)\otimes\dots\otimes G(b_m)\\
+\sum c_{b_1',\cdots,b_m'}^{k_1,\dots,k_m}
G(z^{k_1}b'_1)\otimes\dots\otimes G(z^{k_m}b'_m).
\end{multline*}
Here, the summation ranges over the set of
$(b_1',\dots,b'_m)\in \prod_{\nu=1}^mB_\nu$
and $(k_1,\ldots,k_m)\in\Z^m$
such that
$\sum_{\nu=1}^m k_\nu=0$
and
$k_1+\dots+k_\nu\ge0$ ($\nu=1,\dots,m$).
Moreover we have $c_{b_1',\cdots,b_m'}^{k_1,\dots,k_m}\in \qs\Q[\qs]$.

On the other hand, we have 
\eqn
&&\psi(G(z^{k_1}b'_1)\otimes\dots\otimes G(z^{k_m}b'_m))\\
&&\hs{2em}=(a_1^{k_1}\cdots a_m^{k_m})
\psi\bigl(G(b'_1)\otimes\dots\otimes G(b'_m)\bigr)\\
&&\hs{4em}=(a_1/a_2)^{k_1}(a_2/a_3)^{k_1+k_2}\cdots 
\psi\bigl(G(b'_1)\otimes\dots\otimes G(b'_m)\bigr)\in L,
\eneqn
where
$L=\bigoplus_{b_\nu\in B(M_\nu)} 
A\big(G(b_1)\otimes\dots\otimes G(b_m)\big)$.

Hence we have
\[
\psi\Bigl(G\bigl(b_1\otimes\dots \otimes b_m\bigr)\Bigr)\equiv
\psi\bigl(G(b_1)\otimes\dots\otimes G(b_m)\bigr)
\ \mod\, \qs L.\]
Then Nakayama's lemma implies that
$\bigl\{\psi\bigl(G(b_1\otimes\dots \otimes b_m)
\bigr)\set b_\nu\in B(M_\nu)\bigr\}$
generates $(M_1)_{a_1}\otimes \cdots\otimes(M_m)_{a_m}$.
\qed

We can apply the theorem above
to the fundamental representations.

\Theorem
Let $a_\nu\in K$, and $i_\nu\in I_\dz$ {\rm($\nu=1,\ldots,m$)}.
\bi
\item[(i)]
Assume $a_{\nu}/a_{\nu+1}\in A$ for $\nu=1,\dots,m-1$.
Then $W(\varpi_{i_1})_{a_1}\otimes W(\varpi_{i_2})_{a_2}
\otimes \cdots\otimes W(\varpi_{i_m})_{a_m}$
is generated by $u_{\varpi_{i_1}}\otimes\dots\otimes u_{\varpi_{i_m}}$.
\item[(ii)]
Assume $a_{\nu+1}/a_{\nu}\in A$ for $\nu=1,\dots,m-1$.
Then any non-zero 
submodule of 
$W(\varpi_{i_1})_{a_1}\otimes W(\varpi_{i_2})_{a_2}
\otimes \cdots\otimes W(\varpi_{i_m})_{a_m}$
contains $u_{\varpi_{i_1}}\otimes\dots\otimes u_{\varpi_{i_m}}$.
\ei
\end{theorem}

\vspace{10pt}
Since these theorem hold even if we replace
$K$ and $A$ with $\widehat K$ and $\widehat A$,
they imply
the following consequences as shown in \cite{AK}.

\Prop\label{prop:den}
Assume that $M_j$ is a \good\ $\Us$-module
with dominant extremal vector $u_j$.
The normalized $R$-matrix
$$\Rn_{i,j}(x,y)
\colon  (M_i)_x\otimes (M_j)_y\to (M_j)_y
\otimes (M_i)_x
$$ does not have a pole
at $x/y=a\in \widehat A$.
\enprop

Here $\Rn_{i,j}(x,y)$ is the intertwiner
$(M_i)_x\otimes (M_j)_y\to (M_j)_y\otimes (M_i)_x$
so normalized that it sends 
$u_{i}\otimes u_{j}$ to
$u_{j}\otimes u_{i}$.

Let $\psi_{ij}(x,y)$ be the denominator of
$\Rn_{i,j}(x,y)$. Then one has
\eq\label{eq:psiq1}
\psi_{ij}(x,y)\in 1+A[x/y]\qs x/y.
\eneq

For the sake of simplicity,
we assume that $M_j$ as well as its dual $M_j^*$ is a \good\ $\Us$-module,
and let $u_j$ be a dominant extremal vector of $M_j$.

\Prop
\bnum
\item\label{it}
The extremal vector
$u_{1}\otimes\cdots\otimes u_{m}$
generates $(M_1)_{a_1}\otimes\cdots\otimes (M_m)_{a_m}$
if and only if $\Rn_{i,j}(x,y)$ has no pole at
$x/y=a_i/a_j$ for any
$1\le j<i\le m$.
\item
Any non-zero submodule
$(M_1)_{a_1}\otimes\cdots\otimes (M_m)_{a_m}$
contains $u_{1}\otimes\cdots\otimes u_{m}$,
if and only if $\Rn_{i,j}(x,y)$ has no pole at
$x/y=a_i/a_j$ for any $1\le i<j\le m$.
\item
$(M_1)_{a_1}\otimes\cdots\otimes (M_m)_{a_m}$
is irreducible if and
only if $\Rn_{i,j}(x,y)$
does not have a pole at $x/y=a_i/a_j$ for any 
$1\le i,\,j\le m$ $(i\not =j)$.
\enum
\enprop

\Prop
If $M$ and $M'$ are irreducible finite-dimensional integrable
$\Us$-modules,
then
$M\otimes M'_z$ is an irreducible $\Us$-module 
except for finitely many $z$.
\enprop

\section{Combinatorial $R$-matrices}\label{sec:comb}
Let $M_1$ and $M_2$ be two \good\ $\Us$-modules.
Let $u_\nu$ be the extremal vector of $M_\nu$
with dominant weight ($\nu=1,2$).

Let $\psi(z_1/z_2)$ be the denominator of
the normalized $R$-matrix,
normalized by $\psi\in K[z_1/z_2]$ with $\psi(0)=1$.
Then, by \eqref{eq:psiq1}, we have
\eq\label{eq:psiq}
\psi(z)\in 1+\qs zA[z].
\eneq
We have an intertwiner
$$\psi(z_1/z_2)\Rn\colon (M_1)_\aff\otimes(M_2)_\aff\to
(M_2)_\aff\otimes(M_1)_\aff.$$
We shall first prove the follwoing proposition.

\Prop\label{prop:RnL}
$$\psi(z_1/z_2)\Rn\bigl(L(M_1)_\aff\otimes L(M_2)_\aff\bigr)
\subset L(M_2)_\aff\otimes L(M_1)_\aff.$$
\enprop
\proof
Set $M=(M_1)_\aff\otimes(M_2)_\aff$,
and let $L$ be the smallest crystal lattice
of $M$ containing $A[z_1^{\pm1},z_2^{\pm1}](u_1\otimes u_2)$.
Then $L$ is contained in $L(M)$.
Since every vector in
$B(M)$ is connected with
some $z_1^mu_1\otimes z_2^nu_2$,
$L/\qs L\to L(M)/\qs L(M)$ is surjective.
Hence by the following well-known lemma, there exists 
$g$ such that
\[\text{$g\in 1+\qs A[z_1^{\pm1},z_2^{\pm1}]$ and $gL(M)\subset L$.}
\]
\Lemma
Let $R$ be a commutative ring, $a\in R$
and $F$ a finitely generated $R$-module.
If $F=aF$, then there exists $b\in 1+aR$ such that 
$bF=0$.
\enlemma
Let us define $M'$ and $L'$ in the similar way as $M$ and $L$
by exchanging $M_1$ and $M_2$.
The operator
$T=\psi(z_1/z_2)\Rn\colon M\to M'$ 
commutes with $\te_i$, $\tf_i$, $z_1$, $z_2$,
and it satisfies $T(u_1\otimes u_2)\in L(M')$  by
\eqref{eq:psiq}.
Hence we have
$$TL\subset L(M').$$
Taking $g$ as above, we obtain
$$gT(L(M))\subset TL\subset L(M').$$
Since $L(M')$ is a free $A[z_1^{\pm1},z_2^{\pm1}]$-module
of finite rank, the proposition follows from the following lemma.
\qed

\Lemma
Let $F$ be a free $A[z_1^{\pm1},z_2^{\pm1}]$-module,
and $g$ an element in $1+\qs A[z_1^{\pm1},z_2^{\pm1}]$.
If $\smash{u\in K\mathop\otimes\limits_{A}F}$
satisfies $gu\in F$, then
$u$ belongs to $F$.
\enlemma

Since the proof is elementary, we do not give its proof.

As a corollary of Proposition~\ref{prop:RnL} and \eqref{eq:psiq},
we obtain

\Cor\label{cor:RL}
$$\Rn\Bigl(L\bigl((M_1)_\aff\ct (M_2)_\aff\bigr)\Bigr)
\subset L\bigl((M_2)_\aff\rct L(M_1)_\aff\bigr).$$
\encor

\Conj
$$\psi(z_1/z_2)\bigl((M_1)_\aff\otimes (M_2)_\aff\bigr)
\subset \U[z_1^{\pm1},z_2^{\pm1}](u_1\otimes u_2).$$
\enconj

\vs{5pt}
Set
\eqn
N&=&\U[z_1^{\pm1},z_2^{\pm1}](u_1\otimes u_2)
\subset (M_1)_\aff\otimes (M_2)_\aff,\\
N'&=&\U[z_1^{\pm1},z_2^{\pm1}](u_2\otimes u_1)
\subset(M_2)_\aff\otimes (M_1)_\aff.
\eneqn
Then $\Rn$ gives an isomorphism
\eqn
&&\Rn\colon N\isoto N'.
\eneqn
In \S~\ref{sec:good}, we saw that $N$ (resp. $N'$) has a crystal base
$(L(N),B(M_1)_\aff\otimes B(M_2)_\aff)$
(resp. $(L(N'),B(M_2)_\aff\otimes B(M_1)_\aff)$).
Hence $\Rn$ induces an isomorphism:
\[\Rc\colon B(M_1)_\aff\otimes B(M_2)_\aff
\isoto B(M_2)_\aff\otimes B(M_1)_\aff.\]
We have 
\eqn
&&\Rc(zb_1\otimes b_2)=(1\otimes z)\Rc(b_1\otimes b_2),\\
&&\Rc(b_1\otimes zb_2)=(z\otimes1)\Rc(b_1\otimes b_2).
\eneqn
Hence we have a commutative diagram:
\eq\label{dia:Rc}
\begin{array}{ccc}
B(M_1)_\aff\otimes B(M_2)_\aff
&\xrightarrow{\Rc}&B(M_2)_\aff\otimes B(M_1)_\aff\\
\big\downarrow&&\big\downarrow\\
B(M_1)\otimes B(M_2)&\xrightarrow
{\ \,\sim\,\ }&B(M_2)\otimes B(M_1).
\end{array}
\eneq

Hence one obtains the following proposition.
\Prop
If $B_1$ and $B_2$ are a crystal
base of a \good\ $\Us$-module,
then
$B_1\otimes B_2\simeq B_2\otimes B_1$.
\enprop
By Corollary~\ref{cor:RL}, we have
\eq
\Rn(G(b_1\otimes b_2))=G(\Rc(b_1\otimes b_2)).
\eneq

Setting $\Rc(b_1\otimes b_2)=b'_2\otimes b'_1$ with $b_\nu$, 
$b'_\nu\in B(M_\nu)_\aff$, we define
\[S(b_1\otimes b_2)=\wt(b'_1)-\wt(b_1)=\wt(b_2)-\wt(b'_2)\in Q.\]
By \eqref{eq:normR}, we have
$S(b_1\otimes b_2)\in Q_+:=\sum_{i\in I}\Z_{\ge0}\alpha_i$.
On the other hand, we have
$S(z_1b_1\otimes b_2)=S(b_1\otimes z_2b_2)=S(b_1\otimes b_2)$, and
hence it induces a map:
\eq\label{eq:S+}
&&S\colon B(M_1)\otimes B(M_2)\to Q_+.
\eneq

This map $S$ is characterized by the following properties
(note that $B(M_1)\otimes B(M_2)$ is connected):
{\allowdisplaybreaks
\eq
&&S(u_1\otimes u_2)=0,\\[7pt]
&&
\begin{array}{l}
S(\tf_i(b_1\otimes b_2))\\
\quad=
\begin{cases}
S(b_1\otimes b_2)+\alpha_i&
\text{if $\tf_i(b_1\otimes b_2)=(\tf_ib_1)\otimes b_2$ and\ }\\ 
&\hfill\text{\phantom{if }$\tf_i(b'_2\otimes b'_1)=(\tf_ib'_2)\otimes
  b'_1$,}
\\[3pt]
S(b_1\otimes b_2)-\alpha_i
&\mbox{if $\tf_i(b_1\otimes b_2)=b_1\otimes(\tf_ib_2)$ }\\
&\hfill\text{\phantom{if }$\tf_i(b'_2\otimes
  b'_1)=b'_2\otimes(\tf_ib'_1)$,}
\\[3pt]
S(b_1\otimes b_2)&\mbox{otherwise.}
\end{cases}\end{array}
\eneq
}

\section{Energy function}\label{sec:energy}
In this section, we assume that $M$ is \good,
and we investigate the properties of $M_\aff{}^{\otimes 2}$.
In this case, we have
\eq\label{eq:Rc}
\Rn=\bio\circ\cjn\colon M_\aff\otimes M_\aff\to M_\aff\rct M_\aff.
\eneq
Here $\bio\colon M_\aff\ct M_\aff\to M_\aff\rct M_\aff$
is given by
\[\bio(v\otimes v')=q^{(\lam,\lam)-(\wt(v),\wt(v'))}\ol{v'}\otimes\ol{v}.\]
Indeed, $\Rn$ and $\bio\circ\cjn$ are
$\U$-linear homomorphisms sending $u\otimes u$ to itself,
$z\otimes1$ to $1\otimes z$ and $1\otimes z$ to $z\otimes1$.
Such a homomorphism is unique.

Similarly the identity being a unique automorphism
of $B(M)^{\otimes2}$,
\eqref{dia:Rc} implies that there exists a unique map 
$H\colon B(M)_\aff{}^{\otimes2}\to\Z$
such that
\[\Rc(b_1\otimes b_2)=(z^{-H(b_1\otimes b_2)}b_1)
\otimes (z^{H(b_1\otimes b_2)}b_2).\]
Hence, one has
\[S(b_1\otimes b_2)=\wt(b_2)-\wt(b_1)+H(b_1\otimes b_2)\delta.\]
We call $H$ the {\em energy function}.
We have
$H((zb_1)\otimes b_2)=H(b_1\otimes (z^{-1}b_2))=H(b_1\otimes b_2)+1$.
It is easy to see that
\[B(M)_\aff{}^{\otimes2}=\sqcup_{n\in\Z}H^{-1}(n)\]
is the decomposition of $B(M)_\aff^{\otimes2}$
into the minimal regular subcrystals invariant by $z\otimes z$
(cf.\ \cite{(KMN)^2}).

Embedding $B(M)$ into $B(M)_\aff$ as in \eqref{emb:aff}, 
the energy function restricted on $B(M)^{\otimes2}$
is also characterized by the following two properties:
\begin{equation*}
\begin{array}{l}
H(v\otimes v)=0\quad\mbox{for any extremal vector $v$ of $B(M)$,}\\[5pt]
H(\tf_i(b_1\otimes b_2))=
\begin{cases}
H(b_1\otimes b_2)&
\mbox{if $i\not=0$ and $\tf_i(b_1\otimes b_2)\not=0$,\qquad}\\
H(b_1\otimes b_2)+1&
\mbox{if $i=0$ and}\\
&\hfill\tf_i(b_1\otimes b_2)=(\tf_ib_1)\otimes b_2\not=0,\\
H(b_1\otimes b_2)-1&
\mbox{if $i=0$ and}\\
&\hfill\tf_i(b_1\otimes b_2)=b_1\otimes(\tf_ib_2)\not=0.
\end{cases}
\end{array}
\end{equation*}

Set
\eqn
&&\ba{rl}
N&:=\U[(z\otimes z)^{\pm1},z\otimes1+1\otimes z]
(u\otimes u)\subset M_\aff^{\otimes 2},\\
N'&:=\Ker\bigl(\Rn-1\colon M_\aff^{\otimes2}\to
K(z\otimes z^{-1})\otimes_{K[z\otimes z^{-1}]}M_\aff^{\otimes2}\bigr).
\ea
\eneqn
Then we have
$N\subset N'\subset M_\aff^{\otimes2}$.
Define $L(N)=L(M_\aff)^{\otimes 2}\cap N$
and similarly for $L(N')$.
Then one has $L(N)/\qs L(N)\subset L(N')/\qs L(N')\subset 
L(M_\aff^{\otimes2})/\qs L(M_\aff^{\otimes2})$.
Set 
$$B_0(M_\aff^{\otimes2}):=
\{b_1\otimes b_2\in B(M_\aff)^{\otimes2};
H(b_1\otimes b_2)=0\}.$$
Then 
\[\{(z^n\otimes 1)b;n\in\Z_{\ge0},\,b\in B_0(M_\aff^{\otimes2})\}
\sqcup
\{(1\otimes z^n)b;n\in\Z_{>0},\,b\in B_0(M_\aff^{\otimes2})\}\]
is a basis of $L(M_\aff^{\otimes2})/\qs L(M_\aff^{\otimes2})$.

Let us define the subset $B'$ of
$L(M_\aff^{\otimes2})/\qs L(M_\aff^{\otimes2})$
by
$$B':=\{(z^n\otimes1+\delta(n\not=0)(1\otimes z^n))b;
n\in\Z_{\ge0},\,b\in B_0(M_\aff^{\otimes2})\}.$$
Here, for a statement $P$, we define $\delta(P)$ by
\eqn
&&\delta(P)=\begin{cases}
1&\text{if $P$ is true,}\\
0&\text{if $P$ is false.}
\end{cases}
\eneqn
Then $B'$ is linearly independent.

\Lemma
We have
$B'\subset L(N)/\qs L(N)$.
Moreover, $(L(N),B')$ and $(L(N'),B')$ are 
a crystal base of $N$ and $N'$, respectively.
\enlemma
\proof
It is enough to show that
$B'\subset L(N)/\qs L(N)$ and
$B'$ is a basis of $L(N')/\qs L(N')$.
Since $B_0(M_\aff^{\otimes2})$ is a minimal subcrystal invariant by
$z\otimes z$,
we have $B_0(M_\aff^{\otimes2})\subset L(N)/\qs L(N)$.
Since $L(N)/\qs L(N)$ is invariant by $z^n\otimes1+1\otimes z^n$,
we have $B'\subset L(N)/\qs L(N)$.
It remains to prove that
$B'$ generates $L(N')/\qs L(N')$.

Since $\Rn=1$ on $N'$, we have
$\Rn=1$ on $L(N')/\qs L(N')$,
and hence
$L(N')/\qs L(N')\subset
F:=\{v\in(L(M_\aff)/\qs L(M_\aff))^{\otimes2}\set\Rn(v)=v\}$.

Since the action of $\Rn$ on 
$(L(M_\aff)/\qs L(M_\aff))^{\otimes2}=\Q^{\oplus B(M_\aff)^{\otimes 2}}$
is given by $\Rc$, we can see easily that
$B'$
is a basis of $F$.
\qed

\Prop\label{prop:NN'}
$N=N'$ and it has a global basis
$\{G(b);b\in B'\}$.
\enprop
\proof
We shall apply Theorem~\ref{th:ex}
for $N$ and $N'$.
Set $$N_\Q:=\U_\Q[(z\otimes z)^{\pm1},z\otimes 1+1\otimes z]
(u\otimes u)\subset M_\aff^{\otimes 2},$$
and similarly for $N'_\Q$.
Set $S=\Wt(M_\aff^{\otimes2})\setminus (W(2\lam)+\Z\delta)$.
For $\xi=2w\lam+n\delta$ ($w\in W$, $n\in\Z$),
setting $H=\oplus_{\nu+\mu=n}
K\,(z^\nu\otimes z^\mu+z^\mu\otimes z^\nu)u_{w\lam}^{\otimes2}$,
we have
$H\subset N_\xi\subset N'_\xi\subset H$.
Hence $N_\xi=N'_\xi=H$,
and the condition (iv) in Theorem~\ref{th:ex}
is satisfied for $N$ and $N'$.
The condition (v) follows from the fact that
the weight of any extremal vector of $B(M_\aff)^{\otimes 2}$
is in $W(2\lam)+\Z\delta$.
Hence all the conditions in Theorem~\ref{th:ex}
are satisfied for $N$ and $N'$,
and both $N$ and $N'$ have a global basis.
These two global bases coincide, and hence $N=N'$.
\qed

\Cor\label{cor:N0}
If $b_1$, $b_2\in B(M)_\aff$ satisfy $H(b_1\otimes b_2)=0$,
then 
$$G(b_1\otimes b_2)\in\U_\Q[(z\otimes z)^{\pm1},z\otimes1+1\otimes
z](u\otimes u).$$
Moreover, denoting by $N_0$ the vector subspace generated by
$\{G({b_1\otimes b_2})\set\linebreak[0]H(b_1\otimes b_2)=0\}$,
one has
\eqn
&&\U_\Q[(z\otimes z)^{\pm1},z\otimes1+1\otimes z](u\otimes u)
=\Q[z\otimes1+1\otimes z]\otimes_\Q N_0,\\
&&M_\aff{}^{\otimes 2}=
\Q[z^{{\pm1}}\otimes1,1\otimes z^{\pm1}]\otimes_{\Q[(z\otimes z)^{\pm1}]}N_0
=\Q[z^{{\pm1}}\otimes1]\otimes_{\Q}N_0.
\eneqn
\encor

\section{Fock space}\label{sec:Fock}
\subsection{Some properties of \good\ modules}\label{ss:good}
In \cite{Fock}, we defined the wedge spaces and the Fock spaces
for a finite-dimensional integrable $\Us$-module $V$.
In that paper, we assumed several conditions on $V$.
In this section, we shall show 
that all those conditions are satisfied whenever $V$ is
a \good\ module with a perfect crystal base.
%
%
In \cite{Fock}, we employed the reversed coproduct.
Adapting the notations to ours,
those conditions read as follows.
We set $N:=\U[(z\otimes z)^{\pm1},z\otimes1+1\otimes z](u\otimes u)\subset
(V_\aff)^{\otimes2}$ with an extremal vector $u$ of $V$
of weight $\lam$.
\begin{itemize}
\item[(G)]$V$ is \good. Let $(L,B)$ be the crystal base of $V$.
\item[(P)] $B$ is a perfect crystal.  
\item[(L)] Let $s\colon Q\to \Z$ be the additive function such that
  $s(\alpha_i)=1$, and $\ell\colon B_\aff\to\Z$ be the function
  defined by
$\ell(b)=s(\wt(b)-\wt(u))$. Then one has
\[H(b_1\otimes b_2)\le0\Rightarrow \ell(b_1)\le \ell(b_2).\]
\item[(D)]$\psi\in 1+\qs z A[z]$.
Here $\psi(x/y)$ is the denominator of the normalized $R$-matrix 
$\Rn\colon V_x\otimes V_y\to V_y\otimes V_x$.
\item[(R)] For every pair $(b_1,b_2)$ in
$B_\aff$
with $H(b_1\otimes b_2)=0$, there exists 
$C_{b_1,b_2}\in N$ of the form
\[C_{b_1,b_2}=G(b_1)\otimes G(b_2)
-\sum_{b'_1,b'_2}a_{b'_1,b'_2}G(b'_1)\otimes G(b'_2).\]
Here the sum ranges over $(b'_1,b'_2)\in B_\aff^2$ such that
\eqn
&&H(b'_1\otimes b'_2)>0,\\
&&\ell(b_1)<\ell(b'_1)\le\ell(b_2),\\
&&\ell(b_1)\le\ell(b'_2)<\ell(b_2),
\eneqn
and the coefficients $a_{b'_1,b'_2}$ belong to $\Q[\qs,\qs^{-1}]$.
\end{itemize}

\begin{theorem}[\cite{Fock}]
We assume {\rm (G), (L), (D)} and {\rm (R)}.
Then the wedge space $\bigwedge\limits^mV_\aff$ has a basis
$\{G(b_1)\wedge\cdots\wedge G(b_m)\}$,
where $(b_1,\ldots,b_m)$ ranges over
$(B_\aff)^m$ with
$H(b_j\otimes b_{j+1})>0$ {\rm ($j=1,\ldots,m-1$)}.
\end{theorem}
For the other consequences and the Fock space, see 
\S~\ref{ss:wedge}, \S~\ref{ss:Fock} and \cite{Fock}.

In this section we shall prove the following theorem.
\begin{theorem}
Assume that $V$ is a \good\ $\Us$-module.
Then all the properties above except {\rm (P)} are satisfied.
\end{theorem}

In fact, we shall prove here a little bit stronger results.
In the sequel, we assume that $V$ is a \good\ $\Us$-module.
The property (D) has already been proved
in \eqref{eq:psiq}.
The following lemma immediately implies (L).

\Lemma
If $H(b_1\otimes b_2)\le0$, then $\wt(b_2)-\wt(b_1)\in Q_+$.
\enlemma
\proof
By \eqref{eq:S+},
we have
$S(b_1\otimes b_2)=\wt(b_2)-\wt(b_1)+H(b_1\otimes b_2)\delta\in Q_+$.
Hence if $H(b_1\otimes b_2)\le0$, then
$\wt(b_2)-\wt(b_1)\in Q_+$.
\qed
In order to prove the remaining property (R),
we shall prove the following result on global bases.

\Prop\label{prop:devG}
Assume $H(b_1\otimes b_2)=0$.
Write
\eq
&&G(b_1\otimes b_2)
=\sum_{b'_1,\,b'_2\in B(M)_\aff}a_{b'_1,b'_2}G(b'_1)\otimes G(b'_2).
\eneq
Then we have
\eqn
&&a_{b_1,b_2}=1,\\
&&a_{b'_2,b'_1}=q^{(\lam,\lam)-(\wt b'_1,\wt b'_2)}\ol{a_{b'_1,b'_2}}.
\eneqn
If $a_{b'_1,b'_2}\not=0$, then
\eqn
&&\wt(b'_1)\in \Bigl(\wt(b_1)+Q_+\Bigr)\cap
\Bigl(\wt(b_2)- Q_+\Bigr),\\
&&\wt(b'_2)\in\Bigl(\wt(b_1)+Q_+\Bigr)\cap
\Bigl(\wt(b_2)-Q_+\Bigr).
\eneqn
Moreover $\wt(b'_1)=\wt(b_1)$ implies $(b'_1,b'_2)=(b_1,b_2)$,
and $\wt(b'_1)=\wt(b_2)$ implies $(b'_1,b'_2)=(b_2,b_1)$
\enprop
\proof
We have seen 
$\wt(b'_1)\in \wt(b_1)+Q_+$,
$\wt(b'_2)\in\wt(b_2)-Q_+$. and
$\wt(b'_1)=\wt(b_1)$ implies $(b'_1,b'_2)=(b_1,b_2)$.

Since $\Rn G(b_1\otimes b_2)=G(b_1\otimes b_2)$,
and $\cjn G(b_1\otimes b_2)=G(b_1\otimes b_2)$,
we have $\bio G(b_1\otimes b_2)=G(b_1\otimes b_2)$
by \eqref{eq:Rc}.
Hence we have
\[G(b_1\otimes b_2)=\sum_{b'_1,b'_2\in B(M)_\aff}
q^{(\lam,\lam)-(\wt b'_1,\wt b'_2)}
\ol{a_{b'_1,b'_2}}G(b'_2)\otimes G(b'_1),\]
which gives 
$a_{b'_2,b'_1}=q^{(\lam,\lam)-(\wt b'_1,\wt b'_2)}\ol{a_{b'_1,b'_2}}$.
Hence we obtain the remaining assertions.
\qed

\Conj
Conjecturally, we have
$H(b'_1\otimes b'_2)\ge0$ if $a_{b'_1,b'_2}\not=0$.
\enconj

\vs{5pt}
Let us set
\eqn
I_+(b)&=&\{b'\in B_\aff;
\wt(b')-\wt(b)\in Q_+\setminus\{0\}\}\sqcup\{b\},\\
I_-(b)&=&\{b'\in B_\aff;
\wt(b)-\wt(b')\in Q_+\setminus\{0\}\}\sqcup\{b\}.
\eneqn
The following lemma immediately implies (R).
\Lemma
 For every pair $(b_1,b_2)$ in
$B_\aff$, there exists 
$C_{b_1,b_2}\in N$ of the form
\[C_{b_1,b_2}=G(b_1)\otimes G(b_2)
-\sum_{b'_1,b'_2}a_{b'_1,b'_2}G(b'_1)\otimes G(b'_2).\]
Here the sum ranges over $(b'_1,b'_2)\in B_\aff^2$ such that
$H(b'_1\otimes b'_2)>0$ and $b_1',b_2'\in I_+(b_1)\cap I_-(b_2)$,
and the coefficients $a_{b'_1,b'_2}$ belong to $\Q[\qs,\qs^{-1}]$.
\enlemma
\proof
We shall prove this by the induction on
$\ell(b_2)-\ell(b_1)$.
Note that the assertion is trivial
when $H(b_1\otimes b_2)>0$.
We may assume $H(b_1\otimes b_2)\le0$.
Then (L) implies $\ell(b_2)-\ell(b_1)\ge0$.

Set $n:=-H(b_1\otimes b_2)$.
Then $H(z^nb_1\otimes b_2)=0$.
Hence $G(z^nb_1\otimes b_2)\in N$.
By Proposition~\ref{prop:devG},
we can write
\eqn
&&G(z^nb_1\otimes b_2)=z^nG(b_1)\otimes G(b_2)
+\sum_{b'_1,b'_2}a_{b'_1,b'_2}G(b'_1)\otimes G(b'_2)
\eneqn
where the sum ranges over
$(b'_1,b'_2)$ with
$b'_1,b'_2\in I_+(z^nb_1)\cap I_-(b_2)$
and $b'_1\not=z^nb_1$.
In particular, one has $\ell(z^{-n}b'_1)>\ell(b_1)$.
Then,
\eqn
&&(z^{-n}\otimes 1+\delta(n>0)1\otimes z^{-n})G(z^nb_1\otimes b_2)\\
&&\quad=G(b_1)\otimes G(b_2)+\delta(n>0)G(z^nb_1)\otimes G(z^{-n}b_2)\\
&&\ \quad+\sum_{(b'_1,b'_2)\in I}a_{b'_1,b'_2}
\bigl(G(z^{-n}b'_1)\otimes G(b'_2)
+\delta(n>0)G(b'_1)\otimes G(z^{-n}b'_2)\bigr)
\eneqn
belongs to $N_\Q=N\cap(M_\aff^{\otimes2})_\Q$.
Hence, modulo $N_\Q$,
$G(b_1)\otimes G(b_2)$
is a linear combination of
$G(z^nb_1)\otimes G(z^{-n}b_2)$ ($n>0$),
$G(z^{-n}b'_1)\otimes G(b'_2)$ and
$G(b'_1)\otimes G(z^{-n}b'_2)$.

When $n>0$, we have
$\ell(z^{-n}b_2)-\ell(z^nb_1)<\ell(b_2)-\ell(b_1)$,
and the induction hypothesis implies that
$G(z^nb_1)\otimes G(z^{-n}b_2)$ is, modulo $N_\Q$, a linear combination of
$G(b''_1)\otimes G(b''_2)$
with
$H(b''_1\otimes b''_2)>0$ and
$b''_1,b''_2\in I_+(z^nb_1)\cap I_-(z^{-n}b_2)
\subset I_+(b_1)\cap I_-(b_2)$.

Similarly, we have
$\ell(b'_2)-\ell(z^{-n}b'_1)<\ell(b_2)-\ell(b_1)$.
Hence, modulo $N_\Q$,
$G(z^{-n}b'_1)\otimes G(b'_2)$ is a linear combination of
$G(b''_1)\otimes G(b''_2)$
with
$H(b''_1\otimes b''_2)>0$ and
$b''_1,b''_2\in I_+(z^{-n}b'_1)\cap I_-(b'_2)
\subset I_+(b_1)\cap I_-(b_2)$.

Finally, since $\ell(z^{-n}b'_2)-\ell(b'_1)
\le\ell(b'_2)-\ell(z^{-n}b'_1)<\ell(b_2)-\ell(b_1)$, the induction
hypothesis implies that
$G(b'_1)\otimes G(z^{-n}b'_2)$ modulo $N_\Q$ is a linear combination of
$G(b''_1)\otimes G(b''_2)$
with
$H(b''_1\otimes b''_2)>0$ and
$b''_1,b''_2\in I_+(b'_1)\cap I_-(z^{-n}b'_2)
\subset I_+(b_1)\cap I_-(b_2)$.
\qed

\renewcommand{\min}{{\qopname \relax m{min}}}
\newcommand{\Gw}{{G^\wedge}}
\newcommand{\symm}{{\on{sym}}}
\newcommand{\F}{{\mathcal{F}}}
\newcommand{\vac}{{{\rm vac}}}
\newcommand{\Gn}{{G^{\,{\rm pure}}}}

\subsection{Wedge spaces}\label{ss:wedge}
Let us recall the construction of the wedge space in \cite{Fock}.
Let $V$ be a \good\ $\Us$-module
with an extremal global basis $u$.
Let us set 
\eq
N&=&\U[(z\otimes z)^{\pm1},z\otimes1+1\otimes z](u\otimes u)
\subset V_\aff^{\otimes2},\\
N_m&=&\sum_{j=0}^{m-2}
V_\aff{}^{\otimes j}\otimes N\otimes V^{\otimes(m-2-j)}
\subset V_\aff^{\otimes m}.
\eneq
The wedge space $\bigwedge^mV_\aff$
is defined by
\eqn
\bigwedge^mV_\aff=V_\aff^{\otimes m}/N_m.
\eneqn
For $v_1,\ldots,v_m\in V_\aff$, let 
$v_1\wedge\cdots\wedge v_m$ denote the image of
$v_1\otimes\cdots\otimes v_m$ 
by the projection $V_\aff^{\otimes m}\to \bigwedge^mV_\aff$.
Let $L(\bigwedge^mV_\aff)\subset 
\bigwedge^mV_\aff$
be the image of $L(V_\aff^{\otimes m})$.
For $b=b_1\otimes\cdots\otimes b_m\in B(V_\aff^{\otimes m})$, we set
$$\Gn(b)=G(b_1)\wedge\cdots\wedge G(b_m).$$

For $b\in B(V_\aff^{\otimes m})$, let $G(b)$ be the global basis of 
$V_\aff^{\otimes m}$ and
$\Gw(b)$ be its image in $\bigwedge^mV_\aff$.
We set
\eqn
&&B(\bigwedge^mV_\aff)=\{b_1\otimes \cdots\otimes b_m
\in B(V_\aff^{\otimes m});\\
&&\hs{90pt}\text{$H(b_\nu\otimes b_{\nu+1})>0$ for $\nu=1,\ldots,m-1$\}.}
\eneqn
Then in \cite{Fock}, the following properties are proved
\bnum
\item
$\{\Gn(b)\set
b\in B(\bigwedge^mV_\aff)\}$
is a basis of $L(\bigwedge^mV_\aff)$.

\item
Identifying
$B(\bigwedge^mV_\aff)$ with a subset of 
$L(\bigwedge^mV_\aff)/\qs L(\bigwedge^mV_\aff)$ by $\Gn$,
$(L(\bigwedge^mV_\aff),B(\bigwedge^mV_\aff))$
is a crystal base of $\bigwedge^mV_\aff$.
\enum

On the other hand, the following proposition follows from
Proposition~\ref{prop:exp:glob}.
\Prop
For $b_1\in B(V_\aff^{\otimes m_1})$ and
$b_2\in B(V_\aff^{\otimes m_2})$, one has
the equality in $V_\aff^{\otimes (m_1+m_2)}$ 
\eqn
&&G(b_1\otimes b_2)
=G(b_1)\otimes G(b_2)+\sum_{b'_1,b'_2}c_{b'_1,b'_2}
G(b'_1)\otimes G(b'_2).
\eneqn
Here the sum ranges over
$(b'_1,b'_2)\in B(V_\aff^{\otimes m_1})\times B(V_\aff^{\otimes m_2})$
such that $\wt(b'_1)-\wt(b_1)=\wt(b_2)-\wt(b'_2)\in
Q_+\setminus\{0\}$,
and the coefficients satisfy $c_{b'_1,b'_2}\in\qs\Q[\qs]$.
\enprop

Set
\eqn
&&
B_0(V_\aff^{\otimes m})
=\{b_1\otimes \cdots\otimes b_m
\in B(V_\aff^{\otimes m});
\text{$H(b_\nu\otimes b_{\nu+1})=0$\quad}\\
&&\hs{200pt}\text{ for $\nu=1,\ldots,m-1\}$,}\\
&&N_m^0=
\bigoplus_{b\in B_0(V_\aff^{\otimes m})}KG(b).
\eneqn

The similar arguments as in Proposition \ref{prop:NN'}
and Corollary \ref{cor:N0} 
show the following proposition.
\Prop\label{prop:gl:wedge}
\eqn
&&\bigl(\U\otimes
\Q[z_1^{\pm1},\ldots,z_m^{\pm1}]^{\,\symm}\bigr)u^{\otimes m}\\
&&\hs{50pt}=\Q[z_1^{\pm1},\ldots,z_m^{\pm1}]^{\symm}
\mathop\otimes\limits_{\Q[(z_1\cdots z_m)^{\pm1}]}
N_m^0.
\eneqn
Here $z_\nu$ is the automorphism of $V_\aff^{\otimes m}$
induced by the action of $z$ on the $\nu$-th factor,
and $\Q[z_1^{\pm1},\ldots,z_m^{\pm1}]^{\,\symm}$ 
is the ring of symmetric Laurent polynomials.

In particular, for any Laurent polynomial $f(z_1,\ldots, z_m)$
symmetric in $(z_\nu, z_{\nu+1})$ for some $\nu$,
$$f(z_1,\ldots, z_m)N_m^0\subset N_m.$$
\enprop

Since $N_m$ is a $\Q[z_1^{\pm1},\ldots,z_m^{\pm1}]^{\,\symm}$-module,
$\bigwedge^mV_\aff$ has a structure of
a $\Q[z_1^{\pm1},\ldots,z_m^{\pm1}]^{\,\symm}$-module.
We denote by $B_n$ the operator on $\bigwedge^mV_\aff$
given by
$\sum_{\nu=1}^mz_\nu^n$.
Then the $B_n$'s commute with one another.

\Lemma\label{lem:ggFock}
For any $b\in B(V_\aff^{\otimes m})$,
one has either $\Gw(b)=0$ or
$\Gw(b)=\pm \Gw(b')$
for some $b'\in B(\bigwedge^mV_\aff)$.
\enlemma
\proof
Set
$b=z^{a_1}b_1\otimes\cdots\otimes z^{a_m}b_m$ with
$H(b_\nu\otimes b_{\nu+1})=0$.
Then we have by Proposition~\ref{prop:gl:wedge}
$$G(b)\equiv
\pm G(z^{a_{\sigma(1)}}b_1\otimes\cdots\otimes z^{a_{\sigma(m)}}b_m)$$
for any permutation $\sigma$.
Hence we may assume that $(a_1,\ldots,,a_m)$ is a decreasing
sequence.
If there is $\nu$ such that $a_\nu=a_{\nu+1}$ then 
$G(b)\in N_m$ by the preceding proposition.
Otherwise $b$ belongs to $B(\bigwedge^mV_\aff)$.
\qed

\subsection{Global basis of the Fock space}\label{ss:Fock}
The purpose of this subsection is to define the global basis of the
Fock space.

Let us now assume that $V$ is a \good\ $\Us$-module with perfect
crystal base $(L,B)$ of level $\ell$.
Let us recall that a simple crystal $B$ is called
{\em perfect} of level $\ell$ if it satisfies the following
conditions.
\bi
\item[(P1)]
Any $b\in B$ satisfies
$\lan c,\eps(b)\ran=\lan c,\vphi(b)\ran\ge\ell$.
Here $\eps(b)=\sum_i\eps_i(b)\cl(\Lambda_i)\in P_\cl$ and
$\vphi(b)=\sum_i\vphi_i(b)\cl(\Lambda_i)\in P_\cl$.
\item[(P2)]
Set $P_\cl^{(\ell)}=\{\lam\in P_\cl\set
\text{$\lan c,\lam\ran=\ell$ and $\lan h_i,\lam\ran\ge0$ for every
  $i$}\}$, the set of dominant weights of level $\ell$, 
and $B_\min=\{b\in B;\lan c,\eps(b)\ran=\ell\}$.
Then  the two maps
$$\eps\colon B_\min\longrightarrow P_\cl^{(\ell)}\quad\text{and}\quad
\vphi\colon B_\min\longrightarrow P_\cl^{(\ell)}$$
are bijective.
\ei
For example, the vector representation of $A^{(1)}_n$
is a \good\ $\Us$-module with a perfect crystal base of level $1$.
Let $(B_\aff)_\min$ be the inverse image of $B_\min$ by the map
$B_\aff\to B$.
Let us take a sequence $\{b_n\}_{n\in\Z}$ in $(B_\aff)_\min$
such that
\eqn
&&\text{$\vphi(b_n\gr)=\eps(b_{n-1}\gr)$ and 
$H(b_n\gr\otimes b_{n-1}\gr)=1$.}
\eneqn
Such a sequence is called a {\em ground state}.
Take a sequence $\{\lam_n\}_{n\in\Z}$ in $P$ such that
$$\text{$\lam_n=\lam_{n-1}+\wt(b_n\gr)$ and 
$\cl(\lam_n)=\vphi(b_n\gr)=\eps(b_{n-1}\gr)$.}$$

In \cite{Fock}, the Fock spaces $\F_r$ ($r\in\Z$)
are constructed, and they satisfy the following properties.
\bfnum
\item
$\F_r$ is an integrable $\U$-module.
\item\label{F2} $\Wt(\F_r)\subset \lam_r+Q_-$.
\item\label{F3} There exist $\Us$-linear endomorphisms
$B_n$ ($n\in\Z\setminus\{0\}$) of $\F_r$
with weight $n\delta$ satisfying the boson commutation
relations
$[B_n,B_m]=\delta_{-n,m}a_n$
for some $a_n\in K\setminus\{0\}$.
\item There exists a $\U$-linear map
$\ccdot\wedge\ccdot\colon\F_r\otimes\bigwedge^mV_\aff\to\F_{r-m}$ such that
$(u\wedge v)\wedge v'=u\wedge (v\wedge v')$ for
$u\in\F_r$, $v\in \bigwedge^mV_\aff$ and
$v'\in \bigwedge^{m'}V_\aff$.
\item
$B_n(u\wedge v)=(B_nu)\wedge v+u\wedge(z^n v)$ for
$n\in \Z\setminus\{0\}$, $u\in \F_r$ and $v\in V_\aff$.
\item
There is a non-zero vector
$\vac_r\in\F_r$ of weight $\lam_r$,
$(\F_r)_{\lam_r}=K\vac_r$.
Moreover one has $\vac_{r+1}\wedge G(b_r\gr)=\vac_r$.
\item\label{F6} $\{u\in\F_r\set \mbox{$B_nu=0$ for any $n>0$ and $e_iu=0$ for
    any $i$}\}$\hb
$\hs{1em}=K\var_r$.
\item \label{F7}
Let $K[B_{-1},B_{-2},\ldots]=K[B_{n};n\not=0]/
\bigl(\sum_{m>0}K[B_{n};n\not=0]B_{m}\bigr)$ be the Fock space of
the boson algebra. Then
$K[B_{-1},B_{-2},\ldots]\otimes V(\lam_r)\isoto \F_r$ as a
$K[B_{n};n\not=0]\otimes\U$-module.
Here $1\otimes u_{\lam_r}$ corresponds to $\vac_r$.
\item
Let $B(\F_r)$ be the set
of sequences $\{b_n\}_{n\ge r}$ satisfying
\eqn
&&\text{$H(b_{n+1}\otimes b_n)>0$ for any $n\ge r$,}\\
&&\text{$b_n=b_n\gr$ for $n\gge r$.}
\eneqn
For $b=\{b_n\}_{n\ge r}\in B(\F_r)$, set
$\Gn(b)=\vac_n\wedge G(b_{n-1})\wedge\cdots\wedge G(b_r)$
for $n\gge r$.
Then
$\{\Gn(b)\set b\in B(\F_r)\}$ is a basis of $\F_r$.
\item
Set $L(\F_r)=\bigcup\limits_{n\ge r}\vac_n\wedge L(\bigwedge^{n-r}V_\aff)$.
Then $(L(\F_r), B(\F_r))$ is a crystal base of $\F_r$.
Here $B(\F_r)$ is identified with a subset of $L(\F_r)/\qs L(\F_r)$ by $\Gn$.
\item
$f_i^{(k)}\vac_r=\vac_{r+1}\wedge G(\tf_i^kb_r\gr)$.
\efnum

\vs{10pt}
Now we shall show 
that the Fock space $\F_r$ has a global basis.

First let us define a bar involution
$\cj$ on $\F_r$ such that
\eq
&&\cj(\vac_r)=\vac_r,\\
&&\text{$[B_n,\cj]=0$ for any $n>0$.}
\eneq
By (F\ref{F7}), there exists a unique
bar involution on $\F_r$
satisfying the conditions above.
Note that
$\cj\circ B_{-n}\circ\cj=\ol{a_n}{a_n}^{-1}B_{-n}$ for $n>0$,
since $[B_n,{a_n}^{-1}B_{-n}]=1$ implies ${a_n}^{-1}B_{-n}$ is $c$-invariant.

We set 
\eqn
&&(\F_r)_\Q=\sum_{m\ge r}\vac_m\wedge \bigwedge^{m-r}(V_\aff)_\Q.
\eneqn

\Lemma\label{lem:Fockgl}
Let $b:=b_1\otimes\cdots\otimes b_m$ be an element of $B_\aff^{\otimes m}$.
\banum
\item
If $H(b_r\gr\otimes b_1)\le0$, then
$\vac_r\wedge \Gw(b)=\vac_r\wedge \Gn(b)=0$ hold in $\F_{r-m}$.
\item
$\vac_{r+1}\wedge \Gw(b_r\gr\otimes b)=\vac_r\wedge\Gw(b)$.
\eanum
\enlemma
\proof
(a) We have
$$G(b)=\sum c_{b'_1,b'}G(b'_1)\otimes G(b'),$$
where the sum ranges over
$b'_1\in B_\aff$ and $b'\in B_\aff{}^{\otimes(m-1)}$
such that $\wt(b'_1)-\wt(b_1)\in Q_+$.
Since $H(b_r\gr\otimes b_1)\le0$,
we have
$\ell(b_{r-1}\gr)<\ell(b_1)\le\ell(b'_1)$
by Lemma 4.2.2 in \cite{Fock}.
Since $\Wt(\F_{r-1})\subset \lam_{r-1}+Q_-$ by (F\ref{F2}), one has
$\vac_r\wedge G(b'_1)=0$. Hence we obtain
$\vac_r\wedge \Gw(b)=0$.
The proof of $\vac_r\wedge \Gn(b)=0$ is similar.
\hb
(b)\ 
The proof is similar.
One has
$$G(b_r\gr\otimes b)=G(b_r\gr)\otimes G(b)
+\sum c_{b'_0,b'}G(b'_0)\otimes G(b'),$$
where the sum ranges over
$b'_0\in B_\aff$ and $b'\in B_\aff{}^{\otimes m}$
such that $\wt(b'_0)-\wt(b_r\gr)\in Q_+\setminus\{0\}$.
Then by the same reasoning on the weight of $\Wt(\F_{r})$,
we have
$\vac_{r+1}\wedge G(b'_0)=0$.
\qed

By the lemma above, for $b=\{b_n\}_{n\ge r}\in B(\F_r)$,
$$G(b):=\vac_m\wedge \Gw(b_{m-1}\otimes\cdots\otimes b_r)$$
does not depend
on $m$ such that $b_j=b_j\gr$ for $j\ge m$.

\Lemma
$\{G(b)\set b\in B(\F_r)\}$ is a basis of the $A$-module $L(\F_r)$.
\enlemma
\proof
Since $b\equiv G(b)\ \mod\,\qs L(\F_r)$,
$\{G(b)\set b\in B(\F_r)\}$ is linearly independent.
Hence it is enough to show that
it generates $L(\F_r)$.

Let $b=(b_1,\ldots,b_m)\in B(\bigwedge^{m} V_\aff)$.
For any integer $N$, we can write
\eqn
&&\Gn(b)=\sum_{b'}a_{b'}\Gw(b')+\sum_{b''}c_{b''}\Gn(b'').
\eneqn
Here $b'$ ranges over $B(V_\aff^{\otimes m})$ and
$b''=b''_1\otimes\cdots\otimes b''_m$ ranges over
$B(V_\aff^{\otimes m})$ with $\ell(b''_1)>N$.
Taking $\ell(b_{m+r-1}\gr)$ as $N$,
one has
$\vac_{m+r}\wedge \Gn(b'')=0$.
Hence one has
$$\vac_{m+r}\wedge \Gn(b)
=\sum_{b'\in B(V_\aff^{\otimes m})}a_{b'}\vac_{m+r}\wedge \Gw(b').$$
Now it is enough to apply Lemma \ref{lem:ggFock}
and Lemma~\ref{lem:Fockgl}.
\qed

\Theorem
$\{G(b)\set b\in B(\F_r)\}$ is a global basis of $\F_r$.
\entheorem
\proof
It remains to prove that the $G(b)$'s are invariant by 
the bar involution $\cj$.
Let $E$ be the vector space over $\Q$ generated by $\{G(b)\set b\in B(\F_r)\}$.
Then $\vac_{r+m}\wedge \Gw(b)$ is contained in $E$
for any $b\in B(V_\aff^{\otimes m})$
by Lemma~\ref{lem:ggFock} and Lemma~\ref{lem:Fockgl}.
We define the involution $\cj'$ of $\F_r$
by
\eqn
&&\text{$\cj'(v)=v$ for any $v\in E$ and}\\
&&\text{$\cj'(av)=\ol{a}\,\cj'(v)$ for any $v\in\F_r$ and $a\in K$.}
\eneqn
We shall show that $\cj'=\cj$.
In order to see this, it is enough to show the following properties:
\eq
&&\cj'(\vac_r)=\vac_r,\label{eq:cj'1}\\
&&\text{$\cj'$ commutes with $B_n$ if $n>0$,}\label{eq:cj'2}\\
&&\text{$\cj'(av)=\ol{a}\cj'(v)$ for any $v\in\F_r$ and $a\in \U$.}
\label{eq:cj'3}
\eneq

The property \eqref{eq:cj'1} is obvious.

Let us first show that $\cj'$ commutes with $B_n$ ($n>0$).
This follows from the fact that
$B_n(\vac_{r+m}\wedge \Gw(b))=\vac_{r+m}\wedge B_n\Gw(b)$
holds for $b\in B(\bigwedge^mV_\aff)$,
and the fact that $B_n\Gw(b)$ belongs to $E$.

Let us show \eqref{eq:cj'3}.
We have evidently $q^h\circ\cj'=\cj'\circ q^{-h}$ for every $h\in P^*$.

The conjugation $\cj'$ commutes with $e_i$,
because, for $b\in B(\F_r)$, $e_iG(b)$ belongs to
$\Q[\qs+\qs^{-1}]\otimes E$.

Finally, let us show that $\cj'$ commutes with $f_i$.
To see this,  we shall prove $f_i\cj'(v)=\cj'(f_iv)$
for any weight vector $v\in\F_r$
by the induction on $\wt(v)$.
For any $j\in I$, one has,
by using the commutativity of $\cj'$ and $e_i$
\eqn
&&e_j(f_i\cj'(v)-\cj'(f_iv))\\
&&\hs{20pt}=(f_ie_j+\delta_{ij}\frac{t_i-t_i^{-1}}{q_i-q_i^{-1}})\cj'(v)
-\cj'\bigl((f_ie_j+\delta_{ij}\frac{t_i-t_i^{-1}}{q_i-q_i^{-1}})v\bigr)\\
&&\hs{40pt}=f_i\cj'(e_jv)-\cj'(f_ie_jv).
\eneqn
Since this vanishes by the induction hypothesis,
$f_i\cj'(v)-\cj'(f_iv)$ is a highest weight vector.
Similarly it is annihilated by all the $B_n$'s ($n>0$).
Since the weight of $f_i\cj'(v)-\cj'(f_iv)$ is not $\lam_r$,
it must vanish by (F\ref{F6}).
Thus we obtain \eqref{eq:cj'3}.
\qed

\Remark
In the case when
$\g=A_n^{(1)}$ and $V$ is the vector representation,
the global basis of the Fock space was introduced by
B. Leclerc and J.-Y. Thibon (\cite{LT,LT2}).
D. Uglov (\cite{U})
generalized this to the case when $\g=A_n^{(1)}\oplus 
A_m^{(1)}$ and $V$ is the tensor product of the vector representations.
The connection of global bases of Fock space and 
Kazhdan-Lusztig polynomials
are also studied by
M. Varagnolo--E. Vasserot (\cite{VV})
and O. Schiffmann (\cite{S}).
\enrem
\section{Conjectural structure of $V(\lam)$}
In this section, we shall present conjectures
that clarify
the structure of $V(\lam)$ and its crystal base $B(\lam)$
for $\lam\in P^0$.
The paper by Beck, Chari and Pressley (\cite{BCP})
should help to solve them.
These conjectures are closely related with those
of G. Lusztig (\cite{LC}).

Let $\lam$ be a dominant integral weight of level $0$.
We write
$\lam=\sum_{i\in I_\dz} m_i\varpi_i$.
Then the module $\otimes_{i\in I_\dz}V(m_i\varpi_i)$
contains the extremal vector
$\mathop\otimes\limits_{i\in I_\dz}u_{m_i\varpi_i}$ whose weight is $\lam$.
Here we can take any ordering of $I_\dz$ to define the tensor product.
Hence we have a
$\U$-linear morphism
$$\Phi_\lam\colon V(\lam)\to \otimes_{i\in I_\dz}V(m_i\varpi_i)$$
sending $u_\lam$ to $\mathop\otimes\limits_{i\in I_\dz}u_{m_i\varpi_i}$.

\Conj
\bnum
\item
$\Phi_\lam$ is a monomorphism.
\item
$\Phi_\lam^{-1}\bigl(\otimes_{i\in I_\dz}L(m_i\varpi_i)\bigr)=L(\lam)$.
\item
By $\Phi_\lam$, we have an isomorphism of crystals
$$B(\lam)\isoto \bigotimes_{i\in I_\dz}B(m_i\varpi_i).$$
\enum
\enconj
Next we shall consider the case when $\lam$ is a multiple of 
a fundamental weight.
There is a morphism of $\U$-modules
$$\Psi_{m,i}\colon V(m\varpi_i)\to V(\varpi_i)^{\otimes m}$$
sending $u_{m\varpi_i}$ to $u_{\varpi_i}^{\otimes m}$.
Let $z_i$ be the $\Us$-linear automorphism of $V(\varpi_i)$
of weight $d_i\delta$ introduced in \S~\ref{ss:fd},
and let $z_\nu$ ($\nu=1,\ldots,m$) be the 
operator of  $V(\varpi_i)^{\otimes m}$ obtained by the action of $z_i$
 on the $\nu$-th factor. It is again
a $\Us$-linear automorphism of $V(\varpi_i)^{\otimes m}$
of weight $d_i\delta$.
Let $B_0(m\varpi_i)$ be the connected component of
$B(m\varpi_i)$ containing $u_{m\varpi_i}$, and
let $B_0(V(\varpi_i)^{\otimes m})$ be the connected component of
$B(\varpi_i)^{\otimes m}$ containing $u_{\varpi_i}^{\otimes m}$.

\Conj
\bnum
\item
$\Psi_{m,i}$ is a monomorphism.
\item
$\Psi_{m,i}^{-1}L(V(\varpi_i)^{\otimes m})=L(m\varpi_i)$.
\item
$B_0(m\varpi_i)\isoto B_0(V(\varpi_i)^{\otimes m})$
by $\Psi_{m,i}$. Moreover the global basis
$G(b)$ with $b\in B_0(m\varpi_i)$ is sent to the corresponding global
basis of $\U u_{\varpi_i}^{\otimes m}
\subset W(\varpi_i)_\aff^{\otimes m}$
constructed in Theorem~\ref{th:glex}.
\item
Let $S$ be the set of Schur Laurent polynomials in $z_1,\ldots,z_m$,
i.e.\ the set of characters of $GL(m)$
($(z_1,\ldots, z_m)$ being the components of the
diagonal matrices).
Then
$\{G(b)\set b\in B(m\varpi_i)\}$
is by $\Psi_{m,i}$ sent to
$\{aG(b);b\in B_0(V(\varpi_i)^{\otimes m}), a\in S\}$.
\enum
\enconj

Note that, for $a$, $a'\in S$ and $b$, 
$b'\in B_0(V(\varpi_i)^{\otimes m})$,
$a\,G(b)=a'\,G(b')$ holds if and only if
$a'=(z_1\cdots z_m)^ra$ and $b=(z_1\cdots z_m)^rb'$ for some $r\in\Z$.

These conjectures imply the following conjecture on $\tU$
analogous to Peter-Weyl theorem.
For $\lam\in P$, let $B_0(\lam)$ be the connected component of $B(\lam)$ 
containing $u_\lam$. Note that if $\lan c,\lam\ran\not=0$, then
$B_0(\lam)=B(\lam)$.
We consider $\bigsqcup_{\lam\in P}B_0(\lam)\times B(-\lam)$
as a crystal over $\g\oplus\g$.
The Weyl group $W$ acts on 
$\bigsqcup_{\lam\in P}B_0(\lam)\times B(-\lam)$
by
$S^*_w\times S^*_w\colon
B_0(\lam)\times B(-\lam)\to B_0(w\lam)\times B(-w\lam)$.

\Conj
$\Bigl(\bigsqcup_{\lam\in P}B_0(\lam)\times B(-\lam)\Bigr)/W
\isoto
B(\tU)$ as a crystal over $\g\times\g$.
\enconj
Here the usual crystal structure on $B(\tU)$
corresponds to the one of $B_0(\lam)$
and the star crystal structure on $B(\tU)$ corresponds to
the one of $B(-\lam)$.
The isomorphism sends $u_\lam\otimes b\in B_0(\lam)\times B(-\lam)$
to $b^*\in B(\tU)$.

\appendix
\newcommand{\bin}[2]{\genfrac{(}{)}{0pt}{}{#1}{#2}}
\newcommand{\hyper}[2]{\genfrac{}{}{0pt}{}{#1}{#2}}

%
%
\section{}\label{app:qhg}
In this appendix, we shall give a proof of
\eqref{eq:anne} due to Anne Schilling.
Let us define
\begin{equation*}
(a)_n = (a;q)_n = \prod_{i=0}^{n-1} (1-a q^i).
\end{equation*}
Then in terms of $(a;q)_n$, the $q$-binomial in this paper
is given as
\begin{equation*}
\qbin{m}{n}=q^{n(n-m)} \frac{(q^2;q^2)_m}{(q^2;q^2)_n(q^2;q^2)_{m-n}}.
\end{equation*}
Hence replacing $n\to 2n$ and $q\to q^{1/2}$ in 
\eqref{eq:anne}, it reads as follows:
\begin{lemma}
\begin{multline}\label{formula}
\sum_{k=0}^m (-1)^k q^{\frac{1}{2}k(k+1-2m)-nm}
\frac{(q^n)_k (q)_{2n+m} (q)_{\ell-m+k}}
{(q)_{m-k} (q)_{2n+k} (q)_k (q)_{\ell-m}}\\
=\sum_{k=0}^m q^{k(\ell-m-n+1)} \frac{(q^n)_k (q^{n+1})_{m-k}}
{(q)_k (q)_{m-k}}.
\end{multline}
\end{lemma}

\proof
Using \cite[I.10]{GR}
\begin{equation*}
(a)_{m-k}=\frac{(a)_m}{(q^{1-m}/a)_k} (-\frac{q}{a})^k q^{\bin{k}{2}-mk},
\end{equation*}
the equation \eqref{formula} may be rewritten in hypergeometric notation as
\begin{multline*}
q^{-nm} \frac{(q^{2n+1})_m}{(q)_m} {}_3 \Phi_2 
\left[\hyper{q^{-m},q^n,q^{\ell-m+1}}{q^{2n+1},0};q \right]\\
= \frac{(q^{n+1})_m}{(q)_m} {}_2\Phi_1\left[
\hyper{q^{-m},q^n}{q^{-m-n}};q^{\ell-m-2n+1}\right].
\end{multline*}
However, this formula readily follows from \cite[III.7]{GR} with
the replacements
\eqn
&&n\to m,\ b\to q^n,\ c\to q^{-n-m},\ z\to q^{\ell-m-2n+1}.
\eneqn
\qed

\newcommand{\bb}{b}
\newcommand{\baa}[3]{{#1}b_1\otimes t_{\lam{#2}}\otimes {#3}b_2}
\newcommand{\baaa}[3]{{#1}b_1\otimes t_{{#2}}\otimes {#3}b_2}

\section{Formulas for the crystal $B(\tU))$}\label{table}

In this table, $b_1\in B(\infty)$, $b_2\in B(-\infty)$, 
$\lam\in P$, $b=b_1\otimes t_\lam\otimes b_2$,
$\lam_i=\lan h_i,\lam\ran$ and $\wt_i(b_1)=\lan h_i,\wt(b_1)\ran$.
\begin{equation*}
\begin{array}{rl}
\bb^*&=
b_1^*\otimes t_{-\lam-\wt(b_1)-\wt(b_2)}\otimes b_2^*\,,
\\[2pt]
\eps_i(\bb)&=\max(\eps_i(b_1),\, \eps_i(b_2)-\lam_i-\wt_i(b_1))\,,
\\[2pt]
\vphi_i(\bb)&=\max(\vphi_i(b_1)+\lam_i+\wt_i(b_2),\, \vphi_i(b_2))\,,
\\[2pt]
\wt^*(\bb)&=\wt(\bb^*)=-\lam_i\,,\\[2pt]
\seps_i(\bb)&=\max(\seps_i(b_1),\, \sphi_i(b_2)+\lam_i)\,,
\\[2pt]
\sphi_i(\bb)&=\max(\seps_i(b_1)-\lam_i,\, \sphi_i(b_2))\,,
\ea
\end{equation*}
\begin{equation*}
\ba{rl}
\te_i\bb&=
\left\{\begin{array}{ll}
{\baa{\te_i}{}{}}&\mbox{if $\vphi_i(b_1)\ge \eps_i(b_2)-\lam_i$\,,}
\\[3pt]
\baa{}{}{\te_i}&\mbox{if $\vphi_i(b_1)< \eps_i(b_2)-\lam_i$\,,}
\end{array}
\right.
\\[6pt]
\tf_i\bb&=
\left\{\begin{array}{ll}
{\baa{\tf_i}{}{}}&\quad\mbox{if $\vphi_i(b_1)>\eps_i(b_2)-\lam_i$\,,}
\\[3pt]
\baa{}{}{\tf_i}&\quad\mbox{if $\vphi_i(b_1)\le\eps_i(b_2)-\lam_i$\,,}
\end{array}\right.
\ea
\end{equation*}
\begin{equation*}
\ba{rl}
\te_i^*\bb&=
\left\{\begin{array}{ll}
{\baa{\te_i^*}{-\alpha_i}{}}&
\ \mbox{if $\seps_i(b_1)\ge \sphi_i(b_2)+\lam_i$\,,}
\\[3pt]
\baa{}{-\alpha_i}{\te_i^*}&
\ \mbox{if $\seps_i(b_1)< \sphi_i(b_2)+\lam_i$\,,}
\end{array}\right.
\\[6pt]
\tf_i^*\bb&=
\left\{\begin{array}{ll}
{\baa{\tf_i^*}{+\alpha_i}{}}&
\ \mbox{if $\seps_i(b_1)>\sphi_i(b_2)+\lam_i$\,,}
\\[3pt]
\baa{}{+\alpha_i}{\tf_i^*}&
\ \mbox{if $\seps_i(b_1)\le\sphi_i(b_2)+\lam_i$\,,}
\end{array}\right.
\ea
\end{equation*}
\begin{equation*}
\ba{rl}
\te_i^\max\bb&=
\baa{\te_i^\max}{}{\te_i{}^c}\\[0pt]
&\hs{5em}\mbox{where $c=\max(\eps_i(b_2)-\vphi_i(b_1)-\lam_i,0)$,}
\\[3pt]
\tf_i^\max\bb&=
\baa{\tf_i{}^{c}}{}{\tf_i^\max}\\[0pt]
&\hs{5em}\mbox{where $c=\max(\vphi_i(b_1)-\eps_i(b_2)+\lam_i,0)$\,,}
\\[6pt]
\te_i^*{}^\max\bb&=
\left\{
\begin{array}{l}
\baa{\te_i^*{}^\max}{-(\sphi_i(b_2)+\lam_i)\alpha_i}%
{\te_i^*{}^{\sphi_i(b_2)-\seps_i(b_1)+\lam_i}}\\[3pt]
\hfill\mbox{if $\seps_i(b_1)-\sphi_i(b_2)-\lam_i\le0$\,,}\\[3pt]
\baa{\te_i^*{}^\max}{-\seps_i(b_1)\alpha_i}{}\\[3pt]
\hfill\mbox{if $\seps_i(b_1)-\sphi_i(b_2)-\lam_i\ge0$\,,}
\end{array}\right.
\\[30pt]
\tf_i^*{}^\max\bb&=
\left\{
\begin{array}{l}
\baa{\tf_i^*{}^{\seps_i(b_1)-\sphi_i(b_2)-\lam_i}}%
{+(\seps_i(b_1)-\lam_i)\alpha_i}%
{\tf_i^*{}^\max}\\[3pt]
\hfill
\mbox{if $\seps_i(b_1)-\sphi_i(b_2)-\lam_i\ge0$\,,}\\[3pt]
\baa{}{+\sphi_i(b_2)\alpha_i}{\tf_i^*{}^\max}\\[3pt]
\hfill\mbox{if $\seps_i(b_1)-\sphi_i(b_2)-\lam_i\le0$\,.}
\end{array}\right.
\end{array}
\end{equation*}
Assume now $b=b_1\otimes t_\lam\otimes u_{-\infty}$. If $b$ is extremal,
\begin{equation*}
\begin{array}{rl}
S_ib&=
\left\{
\begin{array}{l}
\tf_i{}^{\wt_i(b_1)+\lam_i}b_1\otimes
t_\lam\otimes u_{-\infty}
\hfill\mbox{if $\eps_i(b)=0$,}\\[3pt]
\te_i^\max b_1\otimes t_\lam\otimes
\te_i{}^{-\vphi_i(b_1)-\lam_i}u_{-\infty}
\quad\mbox{if $\vphi_i(b)=0$.}
\end{array}\right.
\end{array}\end{equation*}
If $b^*$ is extremal,
\begin{equation*}
\begin{array}{rl}
S^*_ib&=
\left\{
\begin{array}{l}
\tf_i^*{}^{-\lam_i}b_1\otimes t_{s_i\lam}\otimes
u_{-\infty}\hfill\mbox{if $\seps_i(b)=0$,}\\[3pt]
\te_i^*{}^\max b_1\otimes t_{s_i\lam}\otimes
\te_i^*{}^{\lam_i-\seps_i(b_1)}u_{-\infty}
\qquad\mbox{if $\vphi^*_i(b)=0$.}
\end{array}\right.
\end{array}
\end{equation*}

\bibliographystyle{unsrt}
\def\same{\,$\raise3pt\hbox to 25pt{\hrulefill}\,$}

\end{document}